\documentclass{amsart}
\usepackage{amsmath}
\usepackage{amssymb}
\usepackage[all]{xy}
\usepackage{graphicx}
\usepackage{mathrsfs}
\usepackage{mathabx}
\usepackage{mathrsfs}
\usepackage{bbold}

\numberwithin{equation}{section}

\usepackage{url}
\usepackage{longtable,tabu}
\usepackage{float}

\usepackage[latin1]{inputenc}
\usepackage{xspace,amssymb,amsfonts,euscript}
\usepackage{amsthm,amsmath}
\usepackage{palatino}
\usepackage{euscript}
\input xy \xyoption {all}

\usepackage{tikz,environ}
\usetikzlibrary{arrows}
\usetikzlibrary{patterns,snakes}

\RequirePackage{color}
\definecolor{myred}{rgb}{0.75,0,0}
\definecolor{mygreen}{rgb}{0,0.5,0}
\definecolor{myblue}{rgb}{0,0,0.65}

\tikzset{%
  DynNode/.style={circle, inner sep=0pt, draw=black, fill=white,minimum size=20pt},
  Greater/.style={pos=0.65, inner sep=0mm, outer sep=0mm},
  highlight/.style={rectangle,rounded corners,fill=red!15,draw=red,
    fill opacity=0.5,thick},
  bendBelow/.style={bend left=70, looseness=2},
  bendAbove/.style={bend right=70, looseness=2},
  object/.style={circle, fill, inner sep=1.5pt, outer sep=0mm},
  labeling/.style={outer sep=0mm, inner sep=0mm},
  1morph/.style={->, shorten >= 0.5pt, >=stealth'},
  2morph/.style={-implies,double,double equal sign distance,
                 shorten >=2pt, shorten <=3pt},
  spot/.style={color=black, thin, dashed},
  sline/.style={color=blue, line width=2pt},
  tline/.style={color=red, line width=2pt},
  uline/.style={color=green, line width=2pt},
  line/.style={color=#1, line width=2pt},
  line/.default=blue,
  sdot/.style={color=blue, thin, fill},
  tdot/.style={color=red, thin, fill},
  udot/.style={color=green, thin, fill},
  dot/.style={color=#1, thin, fill},
  dot/.default=blue
}

\def\BS{{\EuScript B}}

\def\GS{{\EuScript G}}

\def\PS{{\EuScript P}}

\def\TS{{\EuScript T}}

\def\ESS{{\mathscr E}}
\def\FSS{{\mathscr F}}
\def\GSS{{\mathscr G}}

\def\ISS{{\mathscr I}}

\def\LSS{{\mathscr L}}

\def\PSS{{\mathscr P}}

\def\TSS{{\mathscr T}}

\RequirePackage{ifpdf}
\ifpdf
  \IfFileExists{pdfsync.sty}{\RequirePackage{pdfsync}}{}
  \RequirePackage[pdftex,
   colorlinks = true,
   urlcolor = myblue, 
   citecolor = mygreen, 
   linkcolor = myred, 
   pagebackref,
   bookmarksopen=true]{hyperref}
\else
  \RequirePackage[hypertex]{hyperref}
\fi

\RequirePackage{ae, aecompl, aeguill} 

    \def\CM{{\mathbb{C}}}

    \def\FM{{\mathbb{F}}}
  
  \def\hg{{\mathfrak h}}

    \def\KM{{\mathbb{K}}}

    \def\PM{{\mathbb{P}}}
    \def\QM{{\mathbb{Q}}}
    \def\RM{{\mathbb{R}}}

    \def\ZM{{\mathbb{Z}}}

    \def\AC{{\mathcal{A}}}
    \def\BC{{\mathcal{B}}}
    
    \def\DC{{\mathcal{D}}}

    \def\HC{{\mathcal{H}}}

    \def\MC{{\mathcal{M}}}
    \def\NC{{\mathcal{N}}}
    \def\OC{{\mathcal{O}}}

\def\BS{{\EuScript B}}

\def\GS{{\EuScript G}}

\def\PS{{\EuScript P}}

\def\TS{{\EuScript T}}

\newcommand{\nc}{\newcommand} \newcommand{\renc}{\renewcommand}

\def\a{\alpha}

\def\g{\gamma}

\def\d{\delta}

\renc{\l}{\lambda}

\def\r{\rho}

\newcommand{\rdots}{\mathinner{ \mkern1mu\raise1pt\hbox{.}
    \mkern2mu\raise4pt\hbox{.}
    \mkern2mu\raise7pt\vbox{\kern7pt\hbox{.}}\mkern1mu}}

\def\sgn{{\mathrm{sgn}}}
\def\triv{{\mathrm{triv}}}

\def\reg{{\mathrm{reg}}}

\def\sph{{\mathrm{sph}}}
\def\asph{{\mathrm{asph}}}

\DeclareMathOperator{\Tr}{Tr}

\def\un{\underline}

\def\p{{}^p}

\def\to{\rightarrow}

\def\laongto{\laongrightarrow}

\nc{\triright}{\stackrel{[1]}{\to}}
\nc{\laongtriright}{\stackrel{[1]}{\laongto}}

\nc{\Hb}{H^\bullet}

\nc{\Br}{\mathcal{B}}
\nc{\HotRR}{{}_R\mathcal{K}_R}
\nc{\HotR}{\mathcal{K}_R}
\nc{\excise}[1]{}
\nc{\defect}{\text{df}}
\nc{\h}[1]{\underline{H}_{#1}}

\nc{\Ga}{\mathbb{G}_a} 
\nc{\Gm}{\mathbb{G}_m} 

\nc{\Perv}{{\mathbf{P}}}

\nc{\IH}{{\mathrm{IH}}}

\nc{\ic}{\mathbf{IC}}

\nc{\gl}{{\mathfrak{gl}}}
\renc{\sl}{{\mathfrak{sl}}}
\renc{\sp}{{\mathfrak{sp}}}

\renc{\Im}{\textrm{Im}}

\nc{\HBM}{H^{BM}}

\DeclareMathOperator{\codim}{{\mathrm{codim}}}

\DeclareMathOperator{\Hom}{Hom}
 \DeclareMathOperator{\ch}{ch}

\DeclareMathOperator{\Rep}{\mathrm{Rep}}
\DeclareMathOperator{\Parity}{\mathrm{Par}}

\DeclareMathOperator{\id}{id}

\newtheorem{thm}{Theorem}[section]

\newtheorem{cor}[thm]{Corollary}

\newtheorem{question}[thm]{Question}

\theoremstyle{definition}
\newtheorem{ex}[thm]{Example}

\theoremstyle{remark}
\newtheorem{remark}[thm]{Remark}

\DeclareMathOperator{\Spec}{Spec}

\DeclareMathOperator{\rk}{rank}

\newcommand{\ra}{\rightarrow}

\newcommand{\into}{\hookrightarrow}

\def\pt{{\mathrm{pt}}}

\def\Gr{{\EuScript Gr}}

\nc{\simto}{\stackrel{\sim}{\to}}

\DeclareMathOperator{\Fun}{Fun}

\nc{\ext}{\textrm{ext}}

\nc{\fW}{{}^f W}
\nc{\pdot}{ \bullet_p}
\nc{\wdot}{ \bullet}
\nc{\la}{\langle}
\renc{\ra}{\rangle}

\nc{\Wf}{W}
\nc{\Wa}{\mathcal{W}}
\nc{\Sa}{\mathcal{S}}
\nc{\Wae}{\mathcal{W}^{\mathrm{ext}}}

\nc{\Fr}{\mathrm{Fr}} 
\nc{\Repp}{\Rep_0}
\nc{\Repe}{\mathrm{Rep}_0^{\mathrm{ext}}}

\nc{\AntiS}{\mathrm{AS}}
\nc{\CAS}{\mathcal{AS}}

\nc{\HCe}{\HC^{\mathrm{ext}}}

\nc{\Kar}{\mathrm{Kar}}
\nc{\bs}{\mathrm{BS}}

\nc{\Iw}{\mathrm{Iw}}
\nc{\ev}{\mathrm{ev}}

\newcommand{\bk}{\Bbbk}

\nc{\SL}{\mathrm{SL}}
\nc{\GL}{\mathrm{GL}}
\nc{\Sp}{\mathrm{Sp}}
\nc{\Eeight}{\mathrm{E}_8}
\nc{\Gtwo}{\mathrm{G}_2}

\nc{\He}{\mathrm{H}} 
\nc{\Hee}{\mathrm{H}^\mathrm{ext}} 
\nc{\Zvv}{\mathbb{Z}[v]}
\nc{\Zv}{\mathbb{Z}[v^{\pm 1}]}

\nc{\HCat}{\mathcal{H}} 
\nc{\HCate}{\mathcal{H}^\mathrm{ext}} 

\nc{\pcan}{{}^p\un{h}}

\nc{\Loc}{\mathrm{Loc}}
\nc{\mult}{\mathrm{mult}}
\nc{\geom}{\mathrm{geom}}
\nc{\diag}{\mathrm{diag}}
\newcommand{\ubr}[2]{\underbrace{#1}_{#2}}

\nc{\Dmix}{D^{mix}}

\nc{\ses}{\mathrm{Semis}}
\nc{\Par}{\mathrm{Par}}

\def\Chi{\mathscr{X}}
\nc{\IF}{IF}

\makeatletter
\newcommand\leftdash{\!\rotatebox[origin=c]{-60}{$\dabar@\dabar@\dabar@$}\!}
\newcommand\rightdash{\!\rotatebox[origin=c]{60}{$\dabar@\dabar@\dabar@$}\!}
\makeatother

\def \edgeShift {0.5mm} 
\def \wingLen {0.3cm} 

\title{Parity sheaves and the Hecke category}

\author{Geordie Williamson}

\begin{document}

\maketitle

\begin{center}
 \emph{\small{ Dedicated to Wolfgang Soergel, in admiration. \\
``\dots we're still using your imagination\dots'' (L.~A.~Murray)}}
\end{center}


\section*{Introduction}

One of the first theorems of representation theory is Maschke's
theorem: any representation of a finite group over a
field of characteristic zero is semi-simple. This theorem is ubiquitous
throughout mathematics. (We often use it without realising it; for example,  when we 
write a function of one variable as the sum of an odd and an even
function.)
The next step is Weyl's theorem: any finite-dimensional representation of a
compact Lie group is semi-simple\footnote{Weyl first proved his theorem via
integration over the group to produce an invariant Hermitian form. To
do this he needed the theory of manifolds. One can view his proof as
an early appearance of geometric methods in representation
theory.}. It is likewise fundamental: for the circle group Weyl's theorem is
closely tied to the theory of Fourier series.

Beyond the theorems of Maschke and Weyl lies the realm of
where semi-simplicity fails. Non semi-simple phenomena in 
representation theory were first encountered
when studying the modular (i.e. characteristic $p$)
representations of finite groups. This theory is the next step
beyond the classical theory of the character table, and is important
in understanding the deeper structure of finite groups. A second
example (of fundamental importance throughout
mathematics from number theory to mathematical physics) occurs when
studying the infinite-dimensional representation theory of semi-simple
Lie groups and their $p$-adic counterparts.

Throughout the history of representation theory, geometric methods
have played an important role. Over the last forty years, the theory
of intersection cohomology and perverse sheaves has provided powerful
new tools.
To any complex reductive group is naturally associated several
varieties (e.g. unipotent and nilpotent orbits and their
closures, the flag variety and its Schubert varieties, the affine
Grassmannian and its Schubert varieties \dots). In
contrast to the group itself, these varieties are often
singular. The theory of perverse sheaves provides a 
collection of constructible complexes of sheaves (intersection cohomology
sheaves) on such varieties, and the ``IC data'' associated to
intersection cohomology sheaves (graded dimensions of stalks, total
cohomology, \dots) appears throughout Lie theory.

The first example of the power of this theory is the
Kazhdan-Lusztig conjecture (a theorem of Beilinson-Bernstein and
Brylinski-Kashiwara), which expresses the character of a simple
highest weight module over a complex semi-simple Lie algebra in terms
of IC data of Schubert varieties in the flag variety. This theorem is an important first step towards understanding
the irreducible representations of semi-simple Lie groups. A second example is Lusztig's theory of character sheaves,
which provides a family of conjugation equivariant sheaves
on the group which are fundamental to the study of the characters of
finite groups of Lie type.\footnote{The reader is referred to
Lusztig's contribution \cite{Licm} to these proceedings in 1990 for an
impressive list of applications of IC techniques in representation theory.}

An important aspect of the IC data appearing in representation theory
is that it is computable. 
For example, a key step in the proof of the conjecture of Kazhdan and Lusztig is their
theorem that the IC data attached to Schubert varieties in the flag
variety is encoded in Kazhdan-Lusztig polynomials, which are given by  an explicit combinatorial algorithm involving
only the Weyl group. 
Often this
computability of IC data is thanks to the Decomposition Theorem,
which
asserts the semi-simplicity (with
coefficients of characteristic zero) of a direct image sheaf, and
implies that one can compute IC data via a resolution of singularities.

One can view the appearance of the Decomposition Theorem
throughout representation theory as asserting some form of (perhaps well-hidden)
semi-simplicity. A trivial instance of this philosophy is that
Maschke's theorem is equivalent to the Decomposition Theorem for a
finite morphism. A less trivial example is the tendency of categories
in highest weight representation theory to admit Koszul
gradings; indeed, according to \cite{BGS}, a Koszul ring is ``as
close  to  semisimple  as  a $\ZM$-graded  ring  possibly  can
be''. Since the Kazhdan-Lusztig conjecture and its 
proof, many character formulae have been discovered resembling the
Kazhdan-Lusztig conjecture (e.g. for 
affine Lie algebras, quantum groups at roots of unity, Hecke algebras
at roots of unity, \dots) and these are often accompanied by a Koszul grading.

All of the above character formulae involve representations of objects
defined over $\CM$. On the other hand modular representation theory
has been dominated since 1979 by conjectures (the
Lusztig conjecture \cite{L80} on simple representations of reductive
algebraic groups and the James conjecture \cite{James} on simple
representations of symmetric groups) which would imply that characteristic $p$
representations of algebraic groups and symmetric groups are
controlled by related objects over $\CM$ (quantum groups and
Hecke algebras at a $p^{th}$ root of unity) where character formulae
are given by Kazhdan-Lusztig like formulae.

The Decomposition Theorem fails in general with coefficients in a field of
characteristic $p$, as is already evident from the failure of Maschke's
theorem in characteristic $p$. It was pointed out by Soergel \cite{Soe} (and extended by Fiebig
\cite{F} and Achar-Riche \cite{ARFM}) that, after passage through
 deep equivalences, the Lusztig conjecture is
equivalent to the Decomposition Theorem holding for 
Bott-Samelson resolutions of certain complex Schubert varieties, with
coefficients in a field of characteristic $p$. For a fixed morphism, the 
Decomposition Theorem can only fail in finitely many characteristics, 
which implies that the Lusztig conjecture
holds for large primes\footnote{The first proof of Lusztig's conjecture for $p
  \gg 0$ was obtained as a consequence of works by Kazhdan-Lusztig \cite{KLaffine},
  Lusztig \cite{LuM}, Kashiwara-Tanisaki \cite{KT1} and
  Andersen-Jantzen-Soergel \cite{AJS}.}.  More recently, it was discovered that there
are many large characteristics for which the Decomposition
Theorem fails for Bott-Samelson resolutions \cite{WT}. This led to exponentially large counter-examples to the
expected bounds in the Lusztig conjecture as well as counter-examples to
the James conjecture.

Thus the picture for modular representations is much more complicated
than we thought. Recently it has proven useful (see \cite{Soe,JMW2}) to accept the failure
of the Decomposition Theorem in characteristic $p$ and consider
indecomposable summands of direct image sheaves as interesting objects in their
own right. It was pointed out by Juteau, Mautner and the author that, 
 in examples in representation theory, these summands are
often characterised by simple cohomology parity vanishing
conditions, and are called parity sheaves.

Most questions in representation theory whose answer
involves (or is conjectured to involve) Kazhdan-Lusztig polynomials are controlled by
the Hecke category, a categorification of the Hecke algebra of a
Coxeter system. Thus it seems that the Hecke category is a fundamental
object in representation theory, like a group ring or an enveloping algebra.
The goal of this survey is to provide a motivated
introduction to the Hecke category in both its geometric (via parity
sheaves) and diagrammatic (generators and relations)
incarnations. 

When we consider the Hecke category in characteristic
$p$ it gives rise to an interesting new Kazhdan-Lusztig-like
basis of the Hecke algebra, called the $p$-canonical
basis. The failure of this basis to agree with the Kazhdan-Lusztig
basis measures the failure of the Decomposition Theorem in
characteristic $p$. Conjecturally (and provably in many cases), this basis leads to character formulae for 
simple modules for algebraic groups and symmetric groups which are
valid for all $p$. Its uniform calculation for affine Weyl groups and large $p$ seems 
to me to be one of the most interesting problems in
representation theory.\footnote{See \cite{LWbilliards} for a
  conjecture in a very special
case, which gives some idea of (or at least a lower bound on!) the
complexity of this problem.}

If one sees the appearance of the
Decomposition Theorem and  Koszulity  as some form of
semi-simplicity, then this semi-simplicity fails in many
settings in modular representation theory. However it is tempting to
see the appearance of parity sheaves and the $p$-canonical basis as
a deeper and better hidden layer of semi-simplicity, beyond what we
have previously encountered. Some evidence for this is the fact that some
form of Koszul duality still holds, although here IC sheaves are replaced by parity
sheaves and there are no Koszul rings.
\vspace{-.1cm}
\subsection{Structure of the paper}
\begin{enumerate}
\item In \S \ref{sec:DT+p} we discuss the Decomposition Theorem,
  parity sheaves and the role of intersection forms. We conclude with
  examples of parity sheaves and the failure of the Decomposition
  Theorem with coefficients of characteristic $p$.
\item In \S \ref{sec:HC} we introduce the Hecke category. We explain
  two incarnations of this category (via parity sheaves, and via
  diagrammatics) and discuss its spherical and
  anti-spherical modules. We conclude by
  defining the $p$-canonical basis, giving some examples, and
  discussing several open problems.
\item In \S \ref{sec:KD} we give a bird's eye view of Koszul duality for the Hecke
  category in its classical, monoidal and modular forms.
\end{enumerate}
Although this work is motivated by representation theory, we only touch on applications in remarks. The
reader is referred to \cite{W-Takagi} for a survey of
applications of this material to representation theory.

\subsection*{Acknowledgements} It is a pleasure to thank all my
collaborators as well as H.~Andersen, A.~Beilinson, R.~Bezrukavnikov,
M.~A.~de~Cataldo, J.~Chuang,
M.~Khovanov, L.~Migliorini, P.~Polo,
and R.~Rouquier for all they have taught me, and much
besides. Thanks also to P.~McNamara for permission to include
some calculations from work in progress (in Example
\ref{ex:KS}).




\section{The Decomposition Theorem and Parity Sheaves} \label{sec:DT+p}

The Decomposition Theorem is a beautiful theorem about algebraic
maps. However its statement is technical and it takes some effort to
understand its geometric content. To motivate the Decomposition Theorem and the
definition of parity sheaves, we consider one of the paths that led to its discovery, namely
Deligne's proof of the Weil conjectures \cite{DeligneWeilI}. We must
necessarily be brief;
for more background on the Decomposition Theorem see
\cite{BBD,dCMSurvey,WBourbaki}.

\subsection{Motivation: The Weil conjectures}

Suppose that $X$ is a smooth projective variety defined over a finite field
$\FM_q$.
On $X$ one has the Frobenius endomorphism
$\Fr : X \to X$ and the deepest of the Weil conjectures (``purity'') implies that the eigenvalues of
$\Fr$ on the \'etale cohomology $H^i(X)$\footnote{In this section only
  we will use $H^i(X)$ to denote the \'etale cohomology
group $H^i(X_{\overline{\FM}_q},\overline{\QM}_\ell)$
 of the extension of scalars of $X$ to an algebraic closure
$\overline{\FM}_q$ of $\FM_q$, where $\ell$ is a fixed prime number coprime
to $q$.} are of a very special form
(``Weil numbers of weight $i$''). By the
Grothendieck-Lefschetz trace formula, we have
\[
|X(\FM_{q^m})| = \sum_i (-1)^i\Tr( (\Fr^*)^m : H^i(X) \to H^i(X))
\]
for all $m \ge 1$, where $X(\FM_{q^m})$ denotes the (finite) set of
$\FM_{q^m}$-rational points of $X$. In this way, the Weil conjectures
have remarkable implications for the number of $\FM_{q^m}$-points of
$X$.

How should we go about proving purity? We might relate
the cohomology of $X$ to that of other varieties, slowly expanding the world where the Weil
conjectures hold. A first attempt along these lines might be to consider long exact
sequences associated to open or closed subvarieties of $X$. However
this is problematic because purity no longer holds if
one drops the ``smooth'' or ``proper'' assumption.

For any map $f : X \to Y$ of varieties we have a
push-forward functor $f_*$ and its derived functor $Rf_*$ between
(derived) categories of sheaves on $X$ and $Y$. (In this paper we will
never consider non-derived functors; we will write
$f_*$ instead of $Rf_*$
from now on.) The cohomology of $X$ (with its action of
Frobenius) is computed by $p_* \overline{\QM}_{\ell,X}$, where
$\overline{\QM}_{\ell,X}$ denotes the constant sheaf on $X$ and $p : X
\to \pt$ denotes the projection to a point.

This reinterpretation of what cohomology ``means'' provides a more
promising approach to purity. For any map $f : X \to Y$ we can use the
commutative diagram
\[
\xymatrix{ X \ar[rd]^p \ar[d]_f \\
Y \ar[r]_g & \pt
}
\]
and the isomorphism $p_*\overline{\QM}_{\ell,X} = (g \circ f)_* \overline{\QM}_{\ell,X} =g_* (f_* \overline{\QM}_{\ell,X})$
to factor the calculation of $H^*(X)$ into two steps: we can first
understand $f_* \overline{\QM}_{\ell,X}$; then understand the direct
image of this complex to a point.
One can
think of the complex $f_* \overline{\QM}_{\ell,X}$ as a linearisation of the map $f$. For
example, if $f$ is proper and $y$ is a (geometric) point of $Y$ then the stalk at $y$ is
\[
(f_*\QM_{\ell,X})_y = H^*(f^{-1}(y)).
\]
It turns out that\footnote{after passage to $\overline{\FM}_q$}
$f_* \overline{\QM}_{\ell,X}$ splits as a direct sum of simple
pieces (this is the Decomposition Theorem). Thus, each summand
contributes a piece of the cohomology of $X$, and one can try to
understand them separately. This approach provides the skeleton of
Deligne's proof of the Weil conjectures: after some harmless
modifications to $X$, the theory of Lefschetz pencils provides a
surjective morphism $f : X \to \PM^1$, and one has to show purity for
the cohomology of each of the summands of
$f_*\overline{\QM}_{\ell,X}$ (sheaves on $\PM^1$). Showing the purity
of the cohomology of each summand is the heart of the proof, which we
don't enter into here!

\subsection{The Decomposition Theorem}
We now change setting slightly: from now on we consider complex
algebraic varieties equipped with their classical (metric)
topology and sheaves of $\Bbbk$-vector spaces on them, for some field
of coefficients $\Bbbk$. For such a variety $Y$
and a stratification
\[
Y = \bigsqcup_{\l \in \Lambda} Y_\l
\]
of $Y$ into finitely many locally-closed, smooth and connected subvarieties we denote by
$D^b_\Lambda(Y; \Bbbk)$ the full subcategory of the derived category
of complexes of sheaves of $\Bbbk$-vector spaces
with $\Lambda$-constructible\footnote{i.e. those sheaves whose
  restriction to each $Y_\l$ is a local system}
 cohomology sheaves. We will always assume that our stratification is
 such that $D^b_\Lambda(Y; \Bbbk)$  is preserved under Verdier duality (this is
 the case, for example, if our stratification is given  by the orbits
 of a group). We denote by $D^b_c(Y;\Bbbk)$ the constructible
 derived category: it consists of those complexes which are
 $\Lambda$-constructible for some $\Lambda$ as above. Both
 $D^b_\Lambda(Y; \Bbbk)$  and $D^b_c(Y;\Bbbk)$ are triangulated with
 shift functor $[1]$. For any morphism
 $f : X \to Y$ we have functors
\[
\xymatrix{
D^b_c(X;\Bbbk) \ar@/^/[rr]^{f_*, f_!} && \ar@/^/[ll]^{f^*, f^!}  D^b_c(Y;\Bbbk)
}
\]
satisfying a menagerie of relations (see e.g. \cite{dCMSurvey}).

Consider a proper morphism $f : X \to Y$ of complex algebraic
varieties with $X$ smooth. We consider the constant sheaf $\bk_X$ on $X$ with
values in $\bk$ and its (derived) direct image on $Y$:
\[
f_* \bk_X.
\]
A fundamental problem (which we tried to motivate in the
previous section) is to understand how this complex of sheaves decomposes.
The Decomposition Theorem states that, if $\bk$ is a field of
characteristic zero, then $f_* \bk_X$ is semi-simple in the sense of
perverse sheaves. Roughly speaking, this means that much of the topology
of the fibres of $f$ is
``forced'' by the nature of the singularities of $Y$. More precisely,
if we fix a stratification of $Y$ as above for which $f_* \bk_X$ is
constructible, then we have an isomorphism:
\begin{equation}
  \label{eq:DT}
f_* \bk_X \cong \bigoplus H^*_{\l, \LSS} \otimes_\Bbbk
\ic_\l^{\LSS} .
\end{equation}
Here the (finite) sum is over certain pairs $(\l,\LSS)$ where $\LSS$ is
an irreducible local system on $Y_\l$, $H^*_{\l, \LSS}$ is a graded vector
space, and $\ic_\l^{\LSS}$ denotes IC
extension of $\LSS$. (The complex of sheaves
$\ic_\l^{\LSS}$ is supported on $\overline{Y_\l}$ and extends $\LSS[\dim_\CM
Y_\l]$ in a
``minimal'' way, taking into account singularities. For example, if
$\overline{Y}_\l$ is smooth and $\LSS$ extends to 
a local system $\overline{\LSS}$ on $\overline{Y}_\l$, then
$\ic_\l^{\LSS} = \overline{\LSS}[\dim_\CM Y_\l]$.) If $\LSS = \Bbbk$ is the
trivial local system we sometimes write $\ic_\l$ instead of $\ic_\l^\Bbbk$.

Below we will often consider coefficient fields $\Bbbk$ of positive
characteristic, where in general \eqref{eq:DT}
does not hold. We will say that the Decomposition Theorem \emph{holds}
(resp. \emph{fails}) with $\Bbbk$-coefficients if an isomorphism of the form
\eqref{eq:DT} holds (resp. fails).

\subsection{Parity Sheaves}

Consider $f : X \to Y$, a proper map between complex algebraic varieties, with $X$  smooth. Motivated by
the considerations that led to the Decomposition Theorem we ask:

\begin{question} \label{parityQ}
  Fix a field of coefficients $\Bbbk$.
  \begin{enumerate}
  \item What can one say about the
  indecomposable summands of $f_* \Bbbk_X$?
\item What about the indecomposable summands of $f_* \LSS$, for $\LSS$ a local
  system of $\Bbbk$-vector spaces on $X$?
  \end{enumerate}
\end{question}

(Recall $f_*$ means derived direct image.) If $\Bbbk$ is
of characteristic zero, then (1) has a beautiful
answer: by the Decomposition Theorem, any indecomposable summand is a
shift of an
IC extension of an irreducible local system. The same is true of 
(2) if $\LSS$ is irreducible. (This is
Kashiwara's conjecture, proved by Mochizuki
\cite{MochizukiWild}.)

If the characteristic of $\Bbbk$ is positive this question seems
difficult. However it has a nice answer (in terms of ``parity sheaves'') under restrictions on $X$, $Y$ and $f$.

\begin{remark}
  It seems unlikely that this question will have a good answer as phrased in
  general. It is possible that it does have a good answer
  if one instead works in an appropriate category of motives,
  perhaps with restrictions on allowable maps $f$ and local systems $\LSS$.
\end{remark}

Assume that $Y$ admits a stratification $
Y = \bigsqcup_{\l \in \Lambda} Y_\l$
as above. For $\l \in \Lambda$, let $j_\l : Y_\l \into Y$ denote the inclusion.
A complex $\FSS \in D^b_\Lambda(Y;\Bbbk)$ is \emph{even} if
\begin{equation}
  \label{eq:parcond}
  \HC^i(j_\l^* \FSS) = \HC^i(j_\l^! \FSS) = 0 \quad \text{for $i$ odd,
    and all $\l \in \Lambda.$}
\end{equation}
(Here $\HC^i$ denotes the $i^{th}$ cohomology sheaf of a complex of
sheaves.) A complex $\FSS$ is \emph{odd} if $\FSS[1]$ is even; a
complex is \emph{parity} if it can be written as a sum $\FSS_0 \oplus
\FSS_1$ with $\FSS_0$ (resp. $\FSS_1$) even (resp. odd).

\begin{ex}
The archetypal example of a parity complex is $f_* \Bbbk_X[ \dim_\CM
X]$, where $f$ is proper and $X$ is smooth as above, and $f$ is in
addition \emph{even}: $f_* \Bbbk_X[ \dim_\CM
X]$ is $\Lambda$-constructible and the cohomology of the
fibres of $f$ with $\Bbbk$-coefficients vanishes in odd
degree. (Indeed, in this case, $\Bbbk_X[ \dim_\CM X]$ is Verdier
self-dual, hence so is $f_* \Bbbk_X[ \dim_\CM X]$ (by properness) and
the conditions \eqref{eq:parcond} follow from our assumptions on the
cohomology of the fibres of $f$.)
\end{ex}

We make the following (strong) assumptions on each stratum:
\begin{gather}
\label{eq:assump1}
\text{$Y_\l$ is simply connected;} \\
\label{eq:assump2}
H^i(Y_\l,\Bbbk) = 0 \quad \text{for $i$ odd.}
\end{gather}

\begin{thm} \label{thm:parity_class}
  Suppose that $\FSS$ is indecomposable and parity:
  \begin{enumerate}
  \item The support 
of $\FSS$ is
  irreducible, and hence is equal to $\overline{Y}_\l$ for some $\l \in
  \Lambda$.
\item The restriction of $\FSS$ to $Y_\l$ is isomorphic to a constant
  sheaf, up to a shift.
  \end{enumerate}
Moreover, any two indecomposable parity complexes with equal
  support are isomorphic, up to a shift.
\end{thm}
(The proof of this theorem is not difficult, see \cite[\S 2.2]{JMW2}.)
If $\FSS$ is an indecomposable parity complex with support $Y_\l$ then
there is a unique shift of $\FSS$ making it Verdier self-dual. We
denote it by $\ESS_\l^\Bbbk$ or $\ESS_\l$ and call it a \emph{parity
  sheaf}.

\begin{remark}
The above theorem is a uniqueness statement. In general, there might
be no parity complex with support $\overline{Y}_\l$. A condition guaranteeing existence of a
parity sheaf with support $\overline{Y}_\l$ is that $\overline{Y}_\l$ admit an even
resolution. In all settings we consider
below parity sheaves exist for all strata, and thus are classified in
the same way as IC sheaves.
\end{remark}

\begin{remark}
  In contrast to IC sheaves, parity sheaves are only defined up to
  non-canonical isomorphism.
\end{remark}

Below it will also be important to consider the equivariant
setting. We briefly outline the necessary changes. Suppose that a complex
algebraic group $G$ acts on $Y$ preserving strata. Let
$D^b_{G,\Lambda}(X;\Bbbk)$ denote the $\Lambda$-constructible equivariant derived category
\cite{BLu}. We have the usual menagerie of functors associated to
$G$-equivariant maps $f : X \to Y$ which commute with the ``forget
$G$-equivariance functor'' to $D^b_\Lambda(X;\Bbbk)$. In the
equivariant setting the definition of even, odd and parity objects remain
unchanged. Also Theorem \ref{thm:parity_class} holds, if we require
``equivariantly simply connected'' in \eqref{eq:assump1} and state
\eqref{eq:assump2} with equivariant cohomology.

\subsection{Intersection Forms} In the previous section we saw that,
for any proper even map $f : X \to Y$, the derived direct image $f_* \Bbbk_X$ decomposes into a
direct sum of shifts of parity sheaves. In applications it is
important to know precisely what form this decomposition takes. It
turns out that this is encoded in the ranks of certain intersection
forms associated to the strata of $Y$, as we now explain.

For each stratum $Y_\l$ and point $y \in Y_\l$ we can choose a normal
slice $N$ to the stratum $Y_\l$ through $y$. If we set $F :=
f^{-1}(y)$ and $\widetilde{N} := f^{-1}(N)$ then we have a commutative
diagram with Cartesian squares:
\[
\xymatrix{ F \ar[r] \ar[d] & \widetilde{N} \ar[r] \ar[d] & X \ar[d]^f \\
 \{x \} \ar[r] & N \ar[r] & Y }
\]
Set $d := \dim_\CM \widetilde{N} = \dim_\CM N = \codim_\CM ( Y_\l \subset X)$. The
inclusion $F \into \widetilde{N}$ equips the integral 
homology of $F$ with an intersection form (see \cite[\S 3.1]{JMW2})
\[
\IF_\l^j : H_{d-j}(F,\ZM) \times H_{d+j}(F,\ZM) \to
H_0(\widetilde{N},\ZM) = \ZM \quad \text{for $j \in \ZM$.}
\]

\begin{remark}
  Let us give an intuitive explanation for the intersection form:
  suppose we wish to pair the classes of submanifolds
  of real dimension $d-j$ and $d+j$ respectively. We regard our
  manifolds as sitting in $\widetilde{N}$
  and move them until they are transverse. Because
  $(d-j) + (d+j) = 2d$  (the real dimension of $\widetilde{N}$) they will intersect in a finite number of
  signed points, which we then count to get the result.
\end{remark}
 
  \begin{remark}
    The above intersection form depends
    only on the stratum $Y_\l$ (up to non-unique isomorphism):  given any two points $y, y' \in Y_\l$
    and a (homotopy class of) path from $y$ to $y'$ we get an isometry between the two
    intersection forms.
  \end{remark}

Let us assume that our parity assumptions are in force, and that the
homology $H_*(F,\ZM)$ is free for all $\l$. In this case, for any
field $\Bbbk$ the intersection form over $\Bbbk$ is obtained via
extension of scalars from $\IF_\l^j$. We denote this form by $\IF_\l^j
\otimes_\ZM \Bbbk$. The relevance of these forms to the Decomposition
Theorem is the following:

\begin{thm}[{{\cite[Theorem 3.13]{JMW2}}}] \label{thm:IF}
  The multiplicity of $\ESS_\l[j]$ as a summand of $f_*
  \Bbbk_{X}[\dim_\CM X]$ is equal to the rank of
  $\IF_\l^j \otimes \Bbbk$. Moreover, the Decomposition Theorem holds
  if and only if $\IF_\l^j \otimes \Bbbk$ and $\IF_\l^j \otimes \QM$
  have the same rank, for all $\l \in \Lambda$ and $j \in \ZM$.
\end{thm}

\subsection{Examples} \label{ex}
Our goal in this section is to give some examples
of intersection cohomology sheaves and parity sheaves. Throughout,
$\Bbbk$ denotes our field of coefficients and $p$ denotes its
characteristic. (The strata in some of examples below do not satisfy our parity conditions. In each
example this can be remedied by consideration equivariant sheaves for
an appropriate group action.)

\begin{ex} (A nilpotent cone) \label{ex:nullcone}
  Consider the singular $2$-dimensional quadric cone
\[
X = \{ (x,y,z) \; | \; x^2 = -yz \} \subset \CM^3.
\]
Then $X$ is isomorphic to the cone of nilpotent matrices inside
$\mathfrak{sl}_2(\CM) \cong \CM^3$. Let $0$ denote the unique singular
point of $X$ and $X_\reg = X \setminus \{ 0 \}$ the smooth locus.
%
Consider the blow-up of $X$ at $0$:
\[
f : \widetilde{X} \to X
\]
This is a resolution of singularities which is
isomorphic to the Springer resolution under the above isomorphism of
$X$ with the nilpotent cone. It is an isomorphism over $X_\reg$ and
has fibre $\PM^1$ over $0$. In particular, the stalks of the direct
image of the shifted constant sheaf $f_* \Bbbk_{\widetilde{X}}[2]$ are given by:
\[
  \begin{array}{c|ccc}
 & -2 & -1 & 0 \\ \hline
X_\reg & \Bbbk & 0 & 0 \\
\{ 0 \} & \Bbbk & 0 & \Bbbk
\end{array}
\]
One has an isomorphism of $\widetilde{X}$ with the total space of the
line bundle $\OC(-2)$ on $\PM^1$. Under such an isomorphism, the zero
section corresponds to $f^{-1}(0)$. In particular, $f^{-1}(0)$ has
self-intersection $-2$ inside $\widetilde{X}$. It follows that:
\begin{gather*}
  f_* \Bbbk_{\widetilde{X}}[2] \cong \Bbbk_X[2] \oplus \Bbbk_{\{0\}},
 \text{ if $\Bbbk$ is of characteristic $\ne 2$,}\\
  f_* \Bbbk_{\widetilde{X}}[2]
\text{ is indecomposable, if $\Bbbk$ is of characteristic  $2$.}
\end{gather*}
If $p = 2$, the complex $f_* \Bbbk_{\widetilde{X}}[2]$ is an archetypal example
of a parity sheaf.  For further discussion of this example, see \cite[\S
2.4]{JMW1}.
\end{ex}

\begin{ex}(The first singular Schubert variety) \label{ex:gras}
Let $\Gr_2^4$ denote the Grassmannian of 2-planes inside $\CM^4$. Fix a
two-dimensional subspace $\CM^2 \subset \CM^4$ and let $X \subset \Gr_2^4$
denote the closed subvariety (a Schubert variety)
\[
X = \{ V \in \Gr_2^4 \; | \; \dim(V \cap \CM^2) \ge 1 \}.
\]
It is of dimension 3, with unique singular point $V =
\CM^2 \in \Gr_2^4$. The space
\[
\widetilde{X} = \{ V \in \Gr_2^4, L \subset \CM^2 \; | \; \dim L = 1,
L \subset V \cap \CM^2 \}
\]
is smooth, and the map $f : \widetilde{X} \to X$ forgetting $L$ is a resolution of
singularities. This morphism is ``small'' (i.e. the shifted
direct image sheaf $f_*\Bbbk_{\widetilde{X}}[3]$ coincides with the
intersection cohomology complex for any field $\Bbbk$) and even. Thus in this
example the parity sheaf and intersection cohomology sheaf coincide in
all characteristics.
\end{ex}

\begin{ex} (Contraction of the zero section) Suppose that $Y$ is
  smooth of dimension $> 0$ and that $Y \subset T^*Y$ may be
  contracted to a point (i.e. there exists a map $f : T^*Y \to X$ such that $f(Y) = \{ x \}$
  and $f$ is an isomorphism on the complement of $Y$). In this case $x$ is
  the unique singular point of $X$ and the intersection form at $x$ is
  $- \chi(Y)$, where $\chi(Y)$ denotes the Euler characteristic of
  $Y$. If $Y$ has vanishing odd cohomology then $f$ is even and $f_*
  \Bbbk_{ T^*Y}$ is parity. The Decomposition Theorem holds if and
  only if $p$ does not divide $\chi(Y)$.
 \end{ex}

\begin{ex}(A non-perverse parity sheaf) For $n \ge 1$ consider \label{ex:non-perverse}
 \[
X = \CM^{2n}/(\pm 1) = \Spec \CM[x_ix_j \; | \; 1 \le i, j \le 2n].
\]
If $\widetilde{X}$ denotes the total space of $\OC(-2)$ on
$\PM^{2n-1}$ then we have a
resolution 
\[
f : \widetilde{X} \to X.
\]
It is an isomorphism over $X_\reg = X \setminus \{0 \}$ with fibre
$\PM^{2n-1}$ over $0$. The intersection form
\[
\IF_0^j : H_{2n-j}(\PM^{2n-1},\ZM) \times H_{2n+j}(\PM^{2n-1},\ZM) \to \ZM
\]
is non-trivial only for $j = -2n+2, -2n+4, \dots, 2n-2$ in which case it is
the $1 \times 1$ matrix $(-2)$. Thus $f_* \Bbbk_{\widetilde{X}}$ is indecomposable if $p =
2$. Otherwise we have
\[
f_* \Bbbk_{\widetilde{X}}[2n] \cong \Bbbk_X[2n] \oplus \Bbbk_{0}[2n-2]
\oplus \Bbbk_{0}[2n-4]\dots \oplus \Bbbk_{0}[-2n+2].
\]
Because $f$ is even, $f_* \Bbbk_{\widetilde{X}}[2n]$ is parity. It is
indecomposable (and hence is a parity sheaf) if $p =
2$. The interest of this example is that $f_* \Bbbk_{\widetilde{X}}[2n]$ has
many non-zero perverse cohomology sheaves. (See \cite[\S 3.3]{JMW1} for more on this example.)
\end{ex}

\begin{ex}(The generalised Kashiwara-Saito singularity) \label{ex:KS}
  Fix $d \ge 2$ and consider the variety of linear maps
\[
\begin{array}{ccc}
\begin{array}{c} 
\xymatrix{ \CM^d \ar[r]^A \ar[d]_D & \CM^d \ar[d]^B \\
\CM^d\ar[r]_C & \CM^d } \end{array}
& \begin{array}{c} \text{satisfying}\end{array}
& \begin{array}{c} BA = CD = 0, \\
\rk \left ( \begin{matrix} A \\ D \end{matrix} \right ) \le 1, \\
\rk \left ( \begin{matrix} B & C \end{matrix}\right )  \le 1.
\end{array}
\end{array}
\]
This is a singular variety of dimension $6d - 4$. Let $0 = (0,0,0,0)$
denote its most singular point. Consider
\[
\widetilde{X} = \left \{ (A,B,C,D,H_1,L_2,L_3,L_4) \; \middle |
  \; \begin{array}{c} H_1 \in
\Gr_{d-1}^d, L_i \in \Gr_1^d, \\
A \in \Hom(\CM^d/H_1, L_2), B \in
\Hom(\CM^d/L_2,L_4), \\
C \in \Hom(\CM^d/L_3, L_4),  D \in \Hom(\CM^d /H_1, L_3) \end{array} \right \}
\]
(where $\Gr_i^d$ denotes the Grassmannian if $i$-planes in
$\CM^d$). The natural map 
\[
f : \widetilde{X} \to X
\]
is a resolution. We have $F = f^{-1}(0) \cong
(\PM^{d-1})^4$.  The intersection form
\[
H_{6d-4}(F) \times H_{6d-4}(F) \to \ZM
\]
has elementary divisors $(1,\dots,1,d)$. The
Decomposition Theorem holds if and only if $p \nmid d$.

The $d = 2$ case yields an $8$-dimensional singularity which Kashiwara
and Saito showed is smoothly equivalent to a singularity of a
Schubert variety in the flag variety of $\SL_8$ or a quiver variety of
type $A_5$. It tends to show up as a minimal counterexample to optimistic
hopes in representation theory \cite{KSg,Lc,WilJames,WilRed}. Polo
observed (unpublished) that for any $d$ the above singularities occur in Schubert
varieties for $\SL_{4d}$. This shows that the Decomposition Theorem
can fail for type $A$ Schubert varieties for arbitrarily large $p$.
\end{ex}

\section{The Hecke category} \label{sec:HC}

In this section we introduce the Hecke category, a monoidal category
whose Grothendieck group is the Hecke algebra. If one thinks of the
Hecke algebra as providing Hecke operators which act on
representations or function spaces, then the Hecke category consists
of an extra layer of ``Hecke operators
between Hecke operators''.

\subsection{The Hecke algebra} \label{sec:heckealg}

Let $G$ denote a split reductive group over $\FM_q$, and let $T
\subset B \subset G$ denote a maximal torus and Borel subgroup. For
example, we could take $G = \GL_n$, $B = $ upper triangular
matrices and $T = $ diagonal matrices. The set of $\FM_q$-points,
$G(\FM_q)$, is a finite group (e.g. for $G = \GL_n$, $G(\FM_q)$ is the group
of invertible $n\times n$-matrices with coefficients in
$\FM_q$). Many important finite groups, including ``most'' simple
groups, are close relatives of groups of this form.

A basic object in the representation theory of the finite group $G(\FM_q)$ is the
Hecke algebra
\[
H_{\FM_q} := \Fun_{B(\FM_q) \times B(\FM_q)}( G(\FM_q),\CM)
\]
of complex valued functions on $G(\FM_q)$, invariant under
left and right multiplication by $B(\FM_q)$. This is an algebra under
convolution:
\[
(f * f') (g) := \frac{1} {|B(\FM_q)|} \sum_{h \in G(\FM_q)} f(gh^{-1})f'(h).
\]

\begin{remark} Instead we could replace $G(\FM_q)$ by $G(\KM)$ and
  $B(\FM_q)$ by an Iwahori subgroup of $G(\KM)$ (for a local field $\KM$
  with finite residue field), and obtain the 
affine Hecke algebra (important in the representation theory of
$p$-adic groups).
\end{remark}

Let $W$ denote the Weyl group, $S$ its simple reflections, $\ell : W
\to \ZM_{\ge 0}$ the length function of $W$ with respect to $S$ and
$\le$ the Bruhat order. The Bruhat decomposition
\[
G(\FM_q) = \bigsqcup_{w \in W} B(\FM_q) \cdot w B(\FM_q)
\]
shows that $H_{\FM_q}$ has a basis given by indicator
functions $t_w$ of the subsets $B(\FM_q) \cdot w B(\FM_q)$, for
$w \in W$.

Iwahori \cite{Iwahori} showed that $H_{\FM_q}$ may also be described as the unital
algebra generated by $t_s$ for $s \in S$ subject to the relations
\begin{align*}
  t_s^2 &= (q-1)t_s + q, \\
\ubr{t_st_u \dots}{m_{su} \text{ factors}} &= \ubr{t_ut_s \dots}{m_{su} \text{ factors}}
\end{align*}
where $u \ne s$ in the second relation and $m_{su}$ denotes the order of $su$ in $W$. These relations
depend on $q$ in a uniform way and make sense for any Coxeter group. Thus it makes sense to use these
generators and relations to define a new algebra $H$ over $\ZM[q^{\pm
  1}]$ ($q$ is now a formal variable); thus $H$ specialises to the
Hecke algebra defined above via $q \mapsto |\FM_q|$.

For technical reasons it is useful to adjoin a square root of $q$ 
and regard $H$ as defined over $\ZM[q^{\pm 1/2}]$. We then set $v :=
q^{-1/2}$ and $\delta_s := vt_s$, so that the defining relations of $H$ become
\begin{align*}
  \delta_s^2 &= (v^{-1}-v) \delta_s + 1, \\
\ubr{\delta_s \delta_u \dots}{m_{su} \text{ factors}} &= \ubr{\delta_u \delta_s \dots}{m_{su} \text{ factors}}.
\end{align*}
For any reduced expression $w = st \dots u$ (i.e. any expression for
$w$ using $\ell(w)$ simple reflections) we set $\delta_w := \delta_s \delta_t
\dots \delta_u$. We obtain in this way a well-defined $\ZM[v^{\pm 1}]$-basis $\{ \delta_x \;
| \; x \in W \}$ for $H$, the \emph{standard basis}. (This basis
specialises via $q \mapsto |\FM_q|$ to the indicator functions $t_w$
considered above, up to a power of $v$.)


There is an involution $d : H \to H$ defined via
\[v \mapsto v^{-1} \qquad \text{and} \qquad
\d_s \mapsto \d_s^{-1} = \d_s + (v-v^{-1}). \]
Kazhdan and Lusztig \cite{KLp} (see \cite{SoeKL} for a simple proof) showed that for
all $x \in W$ there exists a unique element $b_x$ satisfying
\begin{align}
  \label{eq:duality} \tag{``self-duality"}
  d(b_x) &= b_x, \\
  \label{eq:deg} \tag{``degree bound"}
  b_x &\in \delta_x + \sum_{y < x} v\ZM[v] \delta_y
\end{align}
where $\le$ is Bruhat order. For example $b_s = \delta_s + v$. The set $\{ b_x \; | \; x \in W \}$
is the \emph{Kazhdan-Lusztig basis} of $H$. The polynomials $h_{y,x}
\in \ZM[v]$ defined via $b_x = \sum h_{y,x} \delta_y$ are
(normalisations of) \emph{Kazhdan-Lusztig polynomials}.

\subsection{The Hecke category: geometric incarnation} \label{hecke:geom}

Grothendieck's function-sheaf correspondence (see e.g. \cite[\S 1]{Laumon})
tells us how we should categorify the Hecke algebra $H_{\FM_q}$. Namely, we should consider an
appropriate category of $B \times B$-equivariant sheaves on $G$, with
the passage to $H_{\FM_q}$ being given by the trace of Frobenius at the rational
points $G(\FM_q)$.\footnote{As is always the case with Grothendieck's
function-sheaf correspondence, this actually categorifies the Hecke
algebras of $G(\FM_{q^m})$ for ``all $m$ at once''.} Below we will use
the fact that the multiplication action of $B$ on $G$ is free, and so instead we
can consider $B$-invariant functions (resp. $B$-equivariant sheaves) on $G/B$.

To avoid technical complications, and to ease subsequent discussion,
we will change setting
slightly.  Let us fix a generalised Cartan matrix $C = (c_{st})_{s,t
  \in S}$ and let $(\hg_{\ZM}, \{ \alpha_s \}_{s \in S}, \{
\alpha^\vee_s \}_{s \in S})$ be a Kac-Moody root datum, so that
$\hg_{\ZM}$ is a free and finitely generated $\ZM$-module, $\alpha_s
\in \Hom(\hg_{\ZM},\ZM)$ are ``roots'' and $\alpha_s^\vee \in \hg_{\ZM}$ are
``coroots'' such that $\langle \alpha_s^\vee, \alpha_t \rangle =
c_{st}$. To this data we may associate a Kac-Moody group $\GS$  (a
group ind-scheme over $\CM$) together with a canonical Borel subgroup
$\BS$ and maximal torus $\TS$. The reader is welcome to take $\GS$ to be a complex reductive
group, as per the following remark. (For applications to
representation theory the case of an affine Kac-Moody
group is important.)

\begin{remark}
  If $\GS$ is a complex reductive group and $\TS \subset \BS \subset \GS$ is a
  maximal torus and Borel subgroup, then we can consider the
  corresponding root datum $(\Chi, R, \Chi^\vee, R^\vee)$ (where
  $\Chi$ denotes the characters of $\TS$, $R$ the roots etc.). If $\{
  \alpha_s \}_{s \in S} \subset R$ denotes the simple roots determined
  by $\BS$ then $(\Chi^\vee, \{\alpha_s \}, \{\alpha^\vee_s \})$ is a
  Kac-Moody root datum. The corresponding Kac-Moody group
  (resp. Borel subgroup and maximal torus) is canonically isomorphic
  to $\GS$ (resp. $\BS$, $\TS$).
\end{remark}

We denote by $\GS/\BS$ the flag variety (a projective variety in
the case of a reductive group, and an ind-projective variety in
general). As earlier, we denote by $W$ the Weyl group, $\ell$ the
length function and $\le$ the Bruhat order. We have the Bruhat
decomposition
\[
\GS/\BS = \bigsqcup_{w \in W} X_w \quad \text{where} \quad X_w :=\BS \cdot w\BS/\BS.
\]
The $X_w$ are isomorphic to affine spaces, and are called \emph{Schubert
cells}. Their closures $\overline{X}_w \subset \GS/\BS$ are
projective (and usually singular), and are called \emph{Schubert varieties}.

Fix a field $\Bbbk$ and consider $D^b_{\BS}(\GS/\BS;\Bbbk)$,
the bounded equivariant derived category with
coefficients in $\Bbbk$ (see e.g. \cite{BLu}).\footnote{By definition, any object of
  $D^b_{\BS}(\GS/\BS;\Bbbk)$ is supported on finitely many Schubert cells,
  and hence has finite-dimensional support.} This a
monoidal category under
convolution: given two complexes $\mathscr{F}, \GSS \in D^b_{\BS}(\GS/\BS;\Bbbk)$
their convolution is
\[
\FSS * \GSS := \mult_* ( \FSS \boxtimes_\BS \GSS ),
\]
where: $\GS \times_\BS \GS/\BS$ denotes the quotient of $\GS \times
\GS/\BS$ by $(gb,g' \BS) \sim (g,bg' \BS)$ for all $g,g' \in \GS$ and
$b \in \BS$;
$\mult : \GS \times_\BS \GS/\BS \to \GS/\BS$ is induced by the multiplication on $\GS$;
and $\FSS \boxtimes_\BS \GS  \in D^b_{\BS}(\GS \times_\BS \GS/\BS;\Bbbk)$ is
obtained via descent from $\FSS \boxtimes \GSS \in
D^b_{\BS^3}(\GS \times \GS/\BS;\Bbbk)$.\footnote{The reader is
  referred to \cite{Springer, Nad} for more detail on this
  construction.} (Note that $\mult$ is proper, and so $\mult_* = \mult_!$.)

\begin{remark}
  If $\GS$ is a reductive group and we work over $\FM_q$ instead of $\CM$, then this definition
  categorifies convolution in the Hecke algebra,
  via the function-sheaf correspondence.
\end{remark}

For any $s \in S$ we can consider the parabolic subgroup
\[
\PS_s := \overline{\BS s\BS} = \BS s\BS \sqcup \BS \subset \GS.
\]
We define the Hecke category (in its geometric incarnation) as follows
\[
\HC_{\geom}^{\Bbbk} := \langle \Bbbk_{\PS_s/\BS} \; | \; s \in S \rangle_{*, \oplus,
  [1], \Kar}.
\]
That is, we consider the full subcategory of
$D^b_{\BS}(\GS/\BS;\Bbbk)$ generated by $\Bbbk_{\PS_s/\BS}$ under convolution ($*$),
direct sums ($\oplus$), homological shifts ($[1]$) and direct summands
($\Kar$, for ``Karoubi'').

\begin{remark}
If we were to work over $\FM_q$, then $\Bbbk_{\PS_s/\BS}$ categorifies the indicator function
  of $\PS_s(\FM_q) \subset \GS(\FM_q)$. The definition of the Hecke category is imitating the fact that the
  Hecke algebra is generated by these indicator functions under
  convolution (as is clear from Iwahori's presentation).
\end{remark}

Let $[ \HC_{\geom}^{\Bbbk} ]_{\oplus}$ denote the split Grothendieck
group\footnote{The split Grothendieck group $[\AC]_{\oplus}$ of an
  additive category is the abelian group generated by symbols
$[A]$ for all $A \in \AC$, modulo the relations $[A] = [A'] + [A'']$
whenever $A \cong A' \oplus A''$.}  of
$\HC_{\geom}^{\Bbbk}$. Because
$\HC_{\geom}^{\Bbbk}$ is a monoidal category, $[\HC_{\geom}^{\Bbbk}]_{\oplus}$ is an
algebra via $[ \FSS ]\cdot [\GSS] = [ \FSS * \GSS]$. We view $[
\HC_{\geom}^{\Bbbk}]_{\oplus}$ as a $\ZM[v^{\pm 1}]$-algebra via $v \cdot [\FSS] :=
[\FSS[1]]$.
Recall the Kazhdan-Lusztig basis element $b_s = \delta_s + v$ for all
$s \in S$ from earlier. The following theorem explains the name ``Hecke category''
and is fundamental to all that follows:

\begin{thm} \label{thm:geom_cat}
  The assignment $b_s \mapsto [\Bbbk_{\PS_s/\BS}[1]]$ for all $s \in S$ yields an isomorphism of
  $\ZM[v^{\pm 1}]$-algebras:
\[
H \simto [ \HC_{\geom}^{\Bbbk}]_{\oplus}.
\]
\end{thm}

(This theorem is easily proved using the theory of parity sheaves, as
will be discussed in the next section.) The inverse to the isomorphism in the theorem is given by the
\emph{character map}
\begin{gather*}
  \ch : [\HC_{\geom}^{\Bbbk}]_{\oplus} \simto H   \qquad
\FSS  \mapsto \sum_{x \in W} \dim_\ZM ( H^*(\FSS_{x\BS/\BS}) ) v^{-\ell(x)} \delta_x
\end{gather*}
where: $\FSS_{x\BS/\BS}$ denotes the stalk of the constructible sheaf
on $\GS/\BS$ at the point $x\BS/\BS$ obtained from $\FSS$ by forgetting
$\BS$-equivariance; $H^*$ denotes cohomology; and $\dim_\ZM H^* := \sum
(\dim H^i) v^{-i} \in   \ZM[v^{\pm 1}]$ denotes graded dimension.

\subsection{Role of the Decomposition Theorem} The category
$D^b_\BS(\GS/\BS;\Bbbk)$, and hence also the Hecke category
$\HC ^{\Bbbk}_\geom$, is an example of a Krull-Schmidt category: every
object admits a decomposition into indecomposable objects; and an
object is indecomposable if and only if its endomorphism ring is
local.

Recall that the objects of $\HC ^{\Bbbk}_\geom$ are the direct
summands of finite direct sums of shifts of objects of the form 
\[
\ESS_{(s,t,\dots,u)} := \Bbbk_{\PS_s/\BS} * \Bbbk_{\PS_t/\BS} * \dots  * \Bbbk_{\PS_u/\BS} \in D^b_{\BS}(\GS/\BS)
\]
for any word $(s, t, \dots, u)$ in $S$. The Krull-Schmidt
property implies that any indecomposable object is isomorphic to a direct
summand of a single $\ESS_{(s,t,\dots,u)}$. Thus in order to
understand the objects of $\HC ^{\Bbbk}_\geom$ it is enough to understand the
summands of $\ESS_{(s,t,\dots,u)}$, for any word as above.

For any such word $(s, t, \dots, u)$ we can consider a
Bott-Samelson space
\[
BS_{(s,t,\dots, u)} := \PS_{s} \times_{\BS} \PS_{t} \times_{\BS} \dots \times_{\BS} \PS_u/\BS
\]
and the (projective) morphism $m : BS_{(s,t,\dots, u)}  \to \GS/\BS$ induced by multiplication.
A straightforward argument (using the proper base change theorem) shows
that we have a canonical isomorphism
\[
\ESS_{(s,t,\dots,u)}  = m_* \Bbbk_{BS_{(s,t,\dots,u)}}.
\]
The upshot:  in order to understand the indecomposable objects in 
$\HC_{\geom}^{\Bbbk}$ it is enough to decompose the 
complexes $m_* \Bbbk_{BS}$, for all expressions $(s,t, \dots, u)$ in
$S$.

\begin{remark}
  If $(s,t,\dots,u)$ is a reduced expression for $w \in W$, then the
  map $m$ provides a resolution of singularities of
  the Schubert variety $\overline{X}_w$. These resolutions are often
  called Bott-Samelson resolutions, which explains our notation.
\end{remark}

If the characteristic of our field is zero then we can appeal to the
Decomposition Theorem to deduce that all indecomposable summands of
$m_* \Bbbk_{BS}$ are shifts of the intersection cohomology complexes
of Schubert varieties. Thus
\begin{equation}
  \label{eq:ics}
\HC_{\geom}^{\Bbbk} = \langle \ic_x \; | \; x \in W \rangle_{\oplus,
  [1]} \quad \text{if $\Bbbk$ is of characteristic 0}  
\end{equation}
where $\ic_x$ denotes the ($\BS$-equivariant) intersection cohomology
sheaf of the Schubert variety $\overline{X}_x$. It is also not
difficult (see e.g. \cite{Springer}) to use \eqref{eq:ics} to deduce
that\footnote{Roughly speaking, the two conditions (``self-duality'' + ``degree bound'') characterising the
Kazhdan-Lusztig basis mirror the two conditions (``self-duality'' +
``stalk vanishing'') characterising the IC sheaf.}
\begin{equation}
  \label{eq:0can}
\ch(\ic_x) = b_x \quad \text{if $\Bbbk$ is of characteristic 0}.
\end{equation}
Thus, when the coefficients are of characteristic zero, the
intersection cohomology sheaves categorify the Kazhdan-Lusztig
basis.

It is known that Bott-Samelson resolutions are even. In particular,
$m_* \Bbbk_{BS}$ is a parity complex. Thus, for arbitrary $\Bbbk$ we can
appeal to Theorem \ref{thm:parity_class} to deduce that all indecomposable summands of
$f_* \Bbbk_{BS}$ are shifts of parity sheaves. Thus
\[
\HC_{\geom}^{\Bbbk} = \langle \ESS_x \; | \; x \in W \rangle_{\oplus,
  [1]}\quad \text{for $\Bbbk$ arbitrary}
\]
where $\ESS_x$ denotes the ($B$-equivariant) parity sheaf of the
Schubert variety $\overline{X}_x$. 

\begin{remark}
Recall that for any map there are only finitely many characteristics
in which the Decomposition Theorem fails. Thus, for a fixed $x$ there will be
  only finitely many characteristics in which $\ESS_x^\Bbbk \ne
  \ic_x^\Bbbk$.
\end{remark}

\subsection{The Hecke category: generators and relations}

The above geometric definition of the Hecke category
is analogous to the original definition of 
the Hecke algebra as an algebra of $B$-biinvariant functions. Now we discuss a description of the Hecke category via
generators and relations;  this description is analogous to the
Iwahori presentation of the Hecke algebra.\footnote{The Iwahori
  presentation can be given on two lines. Unfortunately all current
  presentations of the Hecke category need more than two pages!} This
description is due to Elias and the author \cite{EW}, building on work
of Elias-Khovanov \cite{EKh} and Elias \cite{EDC}.

\begin{remark}
  In this section it will be important to keep in mind that monoidal
  categories are fundamentally two dimensional. While group
  presentations (and more generally presentations of categories) 
  occur ``on a line'', presentations of monoidal
  categories (and more generally 2-categories) occur ``in the
  plane''. For background on these ideas the reader is referred to
  e.g. \cite{StreetCat} or \cite[\S 4]{Lauda_sl2}.
\end{remark}

Recall our generalised Cartan matrix $C$, Coxeter system $(W,S$) and
Kac-Moody root datum from earlier. Given $s, t \in S$ we denote by
$m_{st}$ the order (possibly $\infty$) of $st \in W$. We assume:
\begin{equation}
  \label{assump:simply_laced}
  \text{$C$ is simply laced, i.e. $m_{st} \in \{ 2, 3\}$ for $s \ne t$.}
\end{equation}
(We impose this assumption only to shorten the list of relations
below. For the
general case the reader is referred to \cite{EW}.) Recall our ``roots''
and ``coroots''
\[
\{ \a_s \}_{s \in S} \subset \hg_{\ZM}^* \quad \text{and} \quad \{ \a^\vee_s \}_{s \in S} \subset \hg_{\ZM}
\]
such that $\langle \alpha^\vee_s, \alpha_t \rangle = c_{st}$ for all $s,t \in S$. The 
formula $s(v) = v - \langle v, \alpha_s^\vee \rangle \alpha_s$ defines
an action of $W$ on $\hg_{\ZM}^*$. We also assume that our root datum satisfies
that $\alpha_s : \hg_{\ZM} \to \ZM$ and $\alpha_s^\vee : \hg_{\ZM}^* \to
\ZM$ are surjective, for all $s \in S$. (This condition is called
``Demazure surjectivity'' in \cite{EW}. We can always find a Kac-Moody
root datum satisfying this constraint.)

We denote by $R =S(\hg_{\ZM}^*)$ the symmetric algebra of $\hg_{\ZM}^*$ over
$\ZM$. We view $R$ as a graded $\ZM$-algebra with $\deg \hg_{\ZM}^*  =
2$; $W$ acts on $R$ via graded automorphisms. 
For any $s \in S$ we define the \emph{Demazure operator} $\partial_s :
R \to R[-2]$ by
\begin{equation}
  \label{eq:dem}
\partial_s(f) = \frac{f - sf}{\alpha_s}.  
\end{equation}

An \emph{$S$-graph} is a finite, decorated graph, properly embedded
in the planar strip $\mathbb R \times  
[0,1]$, with edges coloured by $S$. The vertices of an $S$-graph are
required to be of the form:
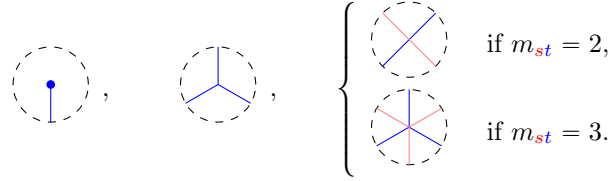
\begin{gather*}
  \begin{array}{c}
\tikz[scale=0.5]{\draw[dashed] (3,0) circle (1cm);
\draw[color=blue] (3,-1) to (3,0);
\node[circle,fill,draw,inner sep=0mm,minimum size=1mm,color=blue] at (3,0) {};}
\end{array},  \qquad 
\begin{array}{c}
\tikz[scale=0.5]{\draw[dashed] (0,0) circle (1cm);
\draw[color=blue] (-30:1cm) -- (0,0) -- (90:1cm);
\draw[color=blue] (-150:1cm) -- (0,0);}
\end{array}, \qquad  
\begin{cases}
\begin{array}{c}\tikz[scale=0.5]{\draw[dashed] (0,0) circle (1cm);
\draw[color=blue] (0,0) -- (45:1cm);
\draw[color=blue] (0,0) -- (-135:1cm);
\draw[color=red!50!white] (0,0) -- (-45:1cm);
\draw[color=red!50!white] (0,0) -- (135:1cm);
}\end{array}
& \text{if $m_{{{\color{red}s}{\color{blue}t}}} = 2,$}\\
\begin{array}{c}\tikz[scale=0.5]{\draw[dashed] (0,0) circle (1cm);
\draw[color=blue] (0,0) -- (90:1cm);
\draw[color=blue] (0,0) -- (-30:1cm);
\draw[color=blue] (0,0) -- (-150:1cm);
\draw[color=red!50!white] (0,0) -- (-90:1cm);
\draw[color=red!50!white] (0,0) -- (30:1cm);
\draw[color=red!50!white] (0,0) -- (150:1cm);
}\end{array} & \text{if $m_{{{\color{red}s}{\color{blue}t}}} = 3.$}\\
\end{cases}
\end{gather*}
The regions (i.e. connected components of the complement of
  our $S$-graph in $\mathbb R \times  
[0,1]$) may be decorated by boxes containing 
  homogeneous elements of $R$.

  \begin{ex}
An $S$-graph (with $m_{{\color{blue}s},
  {\color{red!50!white}t}} = 3$, $m_{{\color{blue}s},
  {\color{green}u}} = 2$, $m_{
  {\color{red!50!white}t}, {\color{green}u}} = 3$, the $f_i \in R$ are homogeneous polynomials):
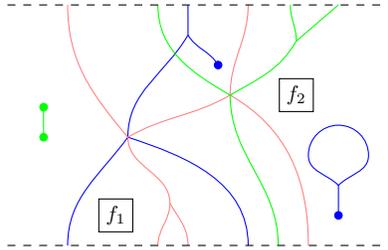
\[
\begin{tikzpicture}[scale=0.8]
  \coordinate (b) at (2.7,0.7);
  \coordinate (a) at (2,1.8);
\coordinate (c) at (3.7,2.5);
\coordinate (dd) at (-0,2);
  \coordinate (e) at (3,3.5);
  \coordinate (ed) at (3.5,3);
\coordinate (f) at (4.8,3.4);
\coordinate (z0) at (0.6,2.3);
\coordinate (z1) at (0.6,1.8);

\node[rectangle,draw,scale=0.8] at (1.8,0.5) {$f_1$};
\node[rectangle,draw,scale=0.8] at (4.8,2.5) {$f_2$};

\draw[blue] (4,0) to[out=90,in=-18] (a);  
\draw[blue] (a) to[out=90,in=-120] (e);
\draw[blue] (e) to (3,4); 
\draw[blue] (ed) to[out=120, in=-60] (e); 
 \node[circle,fill,draw,inner
      sep=0mm,minimum size=1mm,color=blue] at (ed) {};
\draw[blue] (1,0) to [out=90,in=-126] (a); 

\draw[red!50!white] (2.5,0) to[out=90,in=-90] (b) to[out=90,in=-90] (a) to[out=130,in=-90] (1,4);
\draw[red!50!white] (3,0) to[out=90,in=-60] (b);
\draw[red!50!white] (a) to[out=30,in=-150] (c) to[out=90,in=-90] (4,4);
\draw[red!50!white] (5,0) to[out=90,in=-30] (c);

\draw[green] (4.5,0) to[out=90,in=-90] (c) to[out=150,in=-90] (2.5,4);
\draw[green] (c) to[out=30,in=-120] (f) to[out=90,in=-90] (4.7,4);
\draw[green] (f) to (5.5,4);
\draw[green] (z0) to (z1);
 \node[circle,fill,draw,inner
      sep=0mm,minimum size=1mm,color=green] at (z0) {};
 \node[circle,fill,draw,inner
      sep=0mm,minimum size=1mm,color=green] at (z1) {};
\draw[blue] (5.5,1) to[out=135,in=-90] (5,1.5) to[out=90,in=180] (5.5,2) to[out=0,in=90] (6,1.5)
to[out=-90,in=45] (5.5,1) to[out=-90,in=90] (5.5,0.5);
 \node[circle,fill,draw,inner
      sep=0mm,minimum size=1mm,color=blue] at (5.5,0.5) {};

\draw[dashed] (-0,0) -- (6.5,0);
\draw[dashed] (-0,4) -- (6.5,4);
\end{tikzpicture}
\]
\end{ex}

The \emph{degree} of an $S$-graph is the sum over the degrees of
its vertices and boxes, where each box has degree equal to the degree of the
corresponding element of $R$, and the vertices have degrees given by
the following rule: univalent vertices have degree 1, trivalent vertices have degree
$-1$ and $2m_{st}$-valent vertices have degree 0.
 The boundary points of any $S$-graph on $\mathbb R \times\{0\}$
 and on $\mathbb R \times \{1\}$ give\footnote{we read left to right} two words in $S$, called the 
\emph{bottom boundary} and \emph{top boundary}.

\begin{ex}
  The $S$-graph above has degree $0 + \deg f_1 + \deg f_2$. Its bottom
  boundary is $(s,t,t,s,u,s)$ and its top boundary is $(t,u,s,t,u,u)$.
\end{ex}


We are now ready to define a second incarnation of the Hecke category,
which we will denote $\HC_\diag$. By definition, $\HC_\diag$ is monoidally
generated by objects $B_s$, for each $s \in S$. Thus the objects of
$\HC_\diag$ are of the form
\[
B_{(s,t,\dots, u)} := B_sB_t \dots B_u
\]
for some word $(s,t,\dots, u)$ in $S$. (We denote the monoidal
structure in $\HC_{\diag}$ by concatenation.) Thus $\mathbb{1}
:= B_{\emptyset}$ is the monoidal unit.
For any two words $ (s,t, \dots,
u)$ and $(s', t', \dots, v')$ in $S$, $\Hom_{\HC_{\diag}}(B_{(s,t,\dots,u)}, B_{(s',t',\dots,v')})$
is defined to be the free $\ZM$-module
generated by isotopy classes\footnote{i.e. two $S$-graphs are regarded
  as the same
  if one may be obtained from the other by an isotopy of $\RM \times
  [0,1]$ which preserves $\RM \times \{ 0 \}$ and $\RM \times \{ 1 \}$} of $S$-graphs with bottom boundary $ (s,t, \dots,
u)$ and top boundary
$(s', t', \dots, v')$, modulo the local relations below. Composition
(resp. monoidal product) of morphisms is induced by vertical
(resp. horizontal) concatenation of diagrams.

The one colour relations are as follows  (see \eqref{eq:dem} for the definition of $\partial_s$):
\begin{align*}
  \begin{array}{c}
    \tikz[scale=0.6]{\draw[dashed] (0,0) circle (1cm);
\draw[color=red!50!white] (-90:1cm) -- (0,0) -- (90:1cm);
\draw[color=red!50!white] (-0:0.5cm) -- (0,0);
\node[circle,fill,draw,inner sep=0mm,minimum size=1mm,color=red!50!white] at (-0:0.5cm) {};
}
  \end{array}
=
  \begin{array}{c}
    \tikz[scale=0.6]{\draw[dashed] (0,0) circle (1cm);
\draw[color=red!50!white] (-90:1cm) -- (0,0) -- (90:1cm);}
  \end{array} \quad,& \qquad
  \begin{array}{c}
    \tikz[scale=0.6]{\draw[dashed] (0,0) circle (1cm);
\draw[color=red!50!white] (-45:1cm) -- (0.3,0) -- (45:1cm);
\draw[color=red!50!white] (135:1cm) -- (-0.3,0) -- (-135:1cm);
\draw[color=red!50!white] (-0.3,0) -- (0.3,0);
}
  \end{array}
=
  \begin{array}{c}
    \tikz[scale=0.6,rotate=90]{\draw[dashed] (0,0) circle (1cm);
\draw[color=red!50!white] (-45:1cm) -- (0.3,0) -- (45:1cm);
\draw[color=red!50!white] (135:1cm) -- (-0.3,0) -- (-135:1cm);
\draw[color=red!50!white] (-0.3,0) -- (0.3,0);
}
  \end{array}\quad, \\
%
  \begin{array}{c}
    \begin{tikzpicture}[scale=0.6]
      \draw[dashed] (0,0) circle (1cm);
      \draw[red!50!white] (0,0) circle (0.6cm);
      \draw[red!50!white] (0,0.6) --(0,1);
\draw[red!50!white] (0,-0.6) --(0,-1);
    \end{tikzpicture}
  \end{array}= 0 \quad  &,\qquad
 \begin{array}{c}
    \tikz[scale=0.6]{\draw[dashed] (0,0) circle (1cm);
\draw[color=red!50!white] (0,-0.6) -- (0,0.6);
\node[circle,fill,draw,inner sep=0mm,minimum size=1mm,color=red!50!white] at (0,0.6) {};
\node[circle,fill,draw,inner sep=0mm,minimum size=1mm,color=red!50!white] at (0,-0.6) {};
}
  \end{array}
=
  \begin{array}{c}
    \tikz[scale=0.6,rotate=90]{\draw[dashed] (0,0) circle (1cm);
\draw (-0.5,-0.5) rectangle (0.5,0.5);
\node at (0,0) {$\alpha_{\color{red!50!white} s}$};
}
  \end{array} \quad , \\
%
  \begin{array}{c}
    \begin{tikzpicture}[scale=0.6]
      \draw[dashed] (0,0) circle (1cm);
      \draw[red!50!white] (0,1) --(0,-1);
\draw (0.2,-0.4) rectangle (0.7,0.4);
\node at (0.4,0) {$f$};
    \end{tikzpicture}
  \end{array} & = 
  \begin{array}{c}
    \begin{tikzpicture}[scale=0.6]
      \draw[dashed] (0,0) circle (1cm);
      \draw[red!50!white] (0,1) --(0,-1);
\node at (-0.4,0) {${\color{red!50!white} s}f$};
\draw (-0.1,-0.4) rectangle (-0.7,0.4);
    \end{tikzpicture}
  \end{array}
+ 
  \begin{array}{c}
    \begin{tikzpicture}[scale=0.6]
      \draw[dashed] (0,0) circle (1cm);
      \draw[red!50!white] (0,1) --(0,0.6);\node[circle,fill,draw,inner
      sep=0mm,minimum size=1mm,color=red!50!white] at (0,0.6) {}; 
      \draw[red!50!white] (0,-1) --(0,-0.6);\node[circle,fill,draw,inner
      sep=0mm,minimum size=1mm,color=red!50!white] at (0,-0.6) {}; 
\node at (0,0) {$\partial_{\color{red!50!white} s}f$};
\draw (-0.5,-0.4) rectangle (0.5,0.4);
    \end{tikzpicture}
  \end{array} \quad.
\end{align*}

 \begin{remark}
   The first two relations above imply that $B_s$ is a Frobenius
   object in $\HC_\diag$, for all $s \in S$.
 \end{remark}

There are two relations involving two colours. The first is a kind of
``associativity'' (see \cite[(6.12)]{EDC}): 
\[
  \begin{array}{c}
    \begin{tikzpicture}[scale=0.6]
      \draw[dashed] (0,0) circle (1cm);
\draw[red!50!white] (-45:1) -- (135:1);
\draw[blue] (-135:1) -- (45:0.3);
\draw[blue] (45:0.3) -- (20:1);
\draw[blue] (45:0.3) -- (80:1); 
    \end{tikzpicture}
  \end{array}
=
  \begin{array}{c}
    \begin{tikzpicture}[scale=0.6]
      \draw[dashed] (0,0) circle (1cm);
\draw[red!50!white] (-45:1) -- (135:1);
\draw[blue] (-135:1) -- (-135:0.4);
\draw[blue] (-135:0.4) -- (20:1);
\draw[blue] (-135:0.4) -- (80:1); 
    \end{tikzpicture}
  \end{array} \quad \text{if $m_{{{\color{red}s}{\color{blue}t}}} = 2$,} \qquad
 \begin{array}{c}
  \begin{tikzpicture}[scale=0.6]
      \draw[dashed] (0,0) circle (1cm);
\coordinate (a1) at (30:1);\coordinate (a2) at (60:1);\coordinate (a3) at (90:1);
\coordinate (a4) at (135:1);\coordinate (a5) at (-135:1);\coordinate
(a6) at (-90:1); \coordinate (a7) at (-45:1);
\coordinate (l1) at (45:0.4); \coordinate (l0) at (0:0);
\draw[red!50!white] (a5) -- (l0) -- (a3);
\draw[red!50!white] (l0) -- (a7);
\draw[blue] (a1) -- (l1) -- (a2);
\draw[blue] (l1) -- (l0) -- (a4);
\draw[blue] (l0) -- (a6);
  \end{tikzpicture}
\end{array}
= 
\begin{array}{c}
  \begin{tikzpicture}[scale=0.6]
      \draw[dashed] (0,0) circle (1cm);
\coordinate (a1) at (30:1);\coordinate (a2) at (60:1);\coordinate (a3) at (90:1);
\coordinate (a4) at (135:1);\coordinate (a5) at (-135:1);\coordinate
(a6) at (-90:1); \coordinate (a7) at (-45:1);

\coordinate (r1) at (120:0.4); \coordinate (r2) at (-60:0.4); \coordinate (r3) at (-135:0.7);

\draw[blue] (a2) to (r1) to (a4);
\draw[blue] (r1) to[out=-90,in=150] (r2) to (a6);
\draw[blue] (r2) to (a1);

\draw[red!50!white] (a3) to (r1) to[out=-150,in=90] (r3) to (a5);
\draw[red!50!white] (r1) to[out=-30,in=90] (r2) to (a7);
\draw[red!50!white] (r2) to[out=-150,in=0] (r3);

  \end{tikzpicture}
\end{array}\quad \text{if $m_{{{\color{red}s}{\color{blue}t}}} = 3$.}
\]
The second is Elias' ``Jones-Wenzl relation'' (see \cite{EDC}):
\begin{gather*}
  \begin{array}{c}
    \begin{tikzpicture}[scale=0.6]
      \draw[dashed] (0,0) circle (1cm);
\draw[red!50!white] (-45:1) -- (135:1);
\draw[blue] (45:1) -- (-135:0.6);
\node[circle,fill,draw,inner
      sep=0mm,minimum size=1mm,color=blue] at (-135:0.6) {}; 
    \end{tikzpicture}
  \end{array}
=
\begin{array}{c}
    \begin{tikzpicture}[scale=0.6]
      \draw[dashed] (0,0) circle (1cm);
\draw[red!50!white] (-45:1) -- (135:1);
\draw[blue] (45:1) -- (-135:-0.4);
\node[circle,fill,draw,inner
      sep=0mm,minimum size=1mm,color=blue] at (-135:-0.4) {}; 
    \end{tikzpicture}
  \end{array} \qquad \text{if $m_{{{\color{red}s}{\color{blue}t}}} = 2,$} \\
  \begin{array}{c}
    \begin{tikzpicture}[scale=0.6]
      \draw[dashed] (0,0) circle (1cm);
\foreach \r in {30,150,-90}
      \draw[red!50!white] (0,0) -- (\r:1cm);
\foreach \r in {90,-30}
      \draw[blue] (0,0) -- (\r:1cm);
\draw[blue] (0,0) -- (-150:0.7);
\node[circle,fill,draw,inner
      sep=0mm,minimum size=1mm,color=blue] at (-150:0.7) {}; 
    \end{tikzpicture}
  \end{array}
=
  \begin{array}{c}
    \begin{tikzpicture}[scale=0.6]
      \draw[dashed] (0,0) circle (1cm);
\foreach \r in {30,150,-90}
      \draw[red!50!white] (0,0) -- (\r:1cm);
\foreach \r in {90,-30}
{ \draw[blue] (\r:0.5) -- (\r:1cm);
\node[circle,fill,draw,inner sep=0mm,minimum size=1mm,color=blue] at (\r:0.5) {};}
    \end{tikzpicture}
  \end{array} + 
  \begin{array}{c}
    \begin{tikzpicture}[scale=0.6]
      \draw[dashed] (0,0) circle (1cm);
\draw[red!50!white] (-90:1) to[out=90,in=-30] (150:1);
\draw[blue] (90:1) to[out=-90,in=150] (-30:1);
\draw[red!50!white] (30:0.5) -- (30:1cm);
\node[circle,fill,draw,inner sep=0mm,minimum size=1mm,color=red!50!white] at (30:0.5) {};
    \end{tikzpicture}
  \end{array} \qquad \text{if $m_{{{\color{red}s}{\color{blue}t}}} = 3$.}
\end{gather*}

Finally, for each finite standard parabolic subgroup of rank 3 there
is a 3-colour ``Zamolodchikov relation'', which we don't draw here
(see \cite{EW}).
This concludes the definition of $\HC_\diag$. (We remind the reader that if we drop the assumption that $C$ is
simply laced there are more complicated relations, see \cite{EDC,EW}.)
\begin{remark}
    Another way of phrasing the above definition is that $\HC_\diag$ is the
  monoidal category with:
  \begin{enumerate}
  \item generating objects $B_s$ for all $s \in S$;
  \item generating morphisms
  \begin{align*}
      \begin{array}{c}
\tikz[scale=0.5]{
\draw[dashed] (2,-1) to (4,-1); \draw[dashed] (2,1) to (4,1);
\node[rectangle,draw,scale=0.8] at (3,0) {$f$};
} 
\end{array}\in \Hom(\mathbb{1}, \mathbb{1})
\end{align*}
for homogeneous $f \in R$ (recall $\mathbb{1}$ denotes the monoidal unit), as well as
  \begin{align*}
      \begin{array}{c}
\tikz[scale=0.35]{
\draw[dashed] (2,-1) to (4,-1); \draw[dashed] (2,1) to (4,1);
\draw[color=red] (3,-1) to (3,0);
\node[circle,fill,draw,inner sep=0mm,minimum size=1mm,color=red] at (3,0) {};}
\end{array} \in \Hom(B_s, \mathbb{1}), & \qquad
      \begin{array}{c}
\tikz[xscale=0.35,yscale=-0.35]{
\draw[dashed] (2,-1) to (4,-1); \draw[dashed] (2,1) to (4,1);
\draw[color=red] (3,-1) to (3,0);
\node[circle,fill,draw,inner sep=0mm,minimum size=1mm,color=red] at (3,0) {};}
\end{array} \in \Hom(\mathbb{1}, B_s), \\
\begin{array}{c}
\tikz[xscale=0.35,yscale=0.35]{
\draw[dashed] (2,-1) to (4,-1); \draw[dashed] (2,1) to (4,1);
\draw[color=red] (2.3,-1) to[out=90,in=-150] (3,0);
\draw[color=red] (3.7,-1) to[out=90,in=-30] (3,0);
\draw[color=red] (3,0) to (3,1);
}
\end{array} \in \Hom(B_sB_s, B_s), & \qquad
\begin{array}{c}
\tikz[xscale=0.35,yscale=-0.35]{
\draw[dashed] (2,-1) to (4,-1); \draw[dashed] (2,1) to (4,1);
\draw[color=red] (2.3,-1) to[out=90,in=-150] (3,0);
\draw[color=red] (3.7,-1) to[out=90,in=-30] (3,0);
\draw[color=red] (3,0) to (3,1);
}
\end{array} \in \Hom(B_s,B_sB_s)
  \end{align*}
for all $s \in S$ and
  \begin{align*}
      \begin{array}{c}
\tikz[scale=0.35]{
\draw[dashed] (2,-1) to (4,-1); \draw[dashed] (2,1) to (4,1);
\draw[color=red] (2.3,-1) to (3.7,1);
\draw[color=blue] (3.7,-1) to (2.3,1);}
\end{array} \begin{array}{c} \in \Hom(B_sB_t,B_tB_s), \\ \text{ (if $m_{{{\color{red}s}{\color{blue}t}}} = 2$)} \end{array} & \qquad
      \begin{array}{c}
\tikz[scale=0.35]{
\draw[dashed] (2,-1) to (4,-1); \draw[dashed] (2,1) to (4,1);
\draw[color=red] (2.1,-1) to[out=90,in=-150] (3,0) to (3,1); \draw[color=red] (3.9,-1) to[out=90,in=-30] (3,0);
\draw[color=blue] (2.1,1) to[out=-90,in=150] (3,0) to (3,-1); \draw[color=blue] (3.9,1) to[out=-90,in=30] (3,0);}
\end{array} \begin{array}{c} \in \Hom(B_sB_tB_s,B_tB_sB_t), \\ \text{ (if $m_{{{\color{red}s}{\color{blue}t}}} = 3$)} \end{array} 
  \end{align*}
for all pairs $s, t \in S$,
  \end{enumerate}
subject to the above relations (and additional relations encoding isotopy invariance).
\end{remark}

\begin{remark} The above relations are complicated, and perhaps a more
  efficient presentation is possible.  The following is
  perhaps psychologically helpful. Recall that a \emph{standard parabolic subgroup} is a subgroup of $W$
generated by a subset $ I \subset S$, and its
\emph{rank} is $|I|$.
In Iwahori's presentation one has:
  \begin{gather*}
    \text{generators} \leftrightarrow \text{rank 1,} \\
    \text{relations} \leftrightarrow \text{ranks 1, 2.}
  \end{gather*}
In $\HC_\diag$ one has:
  \begin{gather*}
    \text{generating objects} \leftrightarrow \text{rank 1,} \\
    \text{generating morphisms} \leftrightarrow \text{ranks 1, 2,}\\
    \text{relations} \leftrightarrow \text{ranks 1, 2, 3.}
  \end{gather*}
(More precisely, it is only the \emph{finite} standard parabolic
subgroups which contribute at each step.)
\end{remark}

All relations defining $\HC_{\diag}$ are homogeneous for the grading
on $S$-graphs defined above.  Thus $\HC_{\diag}$ is
enriched in graded $\ZM$-modules. We denote by
$\HC^{\oplus,[1]}_{\diag}$ the additive, graded envelope\footnote{Objects are formal sums $F_1[m_1] \oplus F_2[m_2] \oplus \dots \oplus
F_n[m_n]$ where $F_i$ are objects of $\HC_\diag$ and $m_i \in
\ZM$; and morphisms are matrices, determined by the rule that
$\Hom(F[m], F'[m'])$ is the degree $m' - m$ part 
of $\Hom_{\HC_\diag }(F,F')$.}
of $\HC_{\diag}$. Thus $\HC^{\oplus,[1]}_{\diag}$
is an additive category equipped with a ``shift of grading''
equivalence $[1]$, and an isomorphism of graded abelian groups
\[
\Hom_{\HC_\diag}(B,B') = \bigoplus_{m \in \ZM} \Hom_{\HC_\diag^{\oplus,[1]}}(B,B'[m]).
\]
For any field $\Bbbk$, we define
\[
\HC^{\Bbbk,\Kar}_{\diag} := ( \HC^{\oplus,[1]}_{\diag} \otimes_{\ZM} \Bbbk)^\Kar
\]
where $(-)^\Kar$ denotes Karoubi envelope. In other words, $\HC^{\Bbbk,\Kar}_\diag$
is obtained as the additive Karoubi envelope of the extension of
scalars of $\HC_\diag$ to $\Bbbk$.
As for $\HC_\geom$, let us consider the split Grothendieck group
$[\HC_\diag^{\Bbbk, \Kar}]_{\oplus}$ of $\HC_\diag^{\Bbbk, \Kar}$, which we view as a
$\ZM[v^{\pm 1}]$-algebra in the same was as for $\HC^\Bbbk_\geom$
earlier. The following is the analogue of Theorem \ref{thm:geom_cat} in this setting:

\begin{thm}[{{\cite{EW}}}]\label{thm:diag_cat}
  The map $b_s \mapsto [B_s]$ for all $s \in S$ induces an
  isomorphism of $\ZM[v^{\pm 1}]$-algebras:
\[
H \simto [\HC_\diag^{\Bbbk, \Kar}]_{\oplus}.
\]
\end{thm}
The proof is rather complicated diagrammatic algebra,
and involves first producing a basis of morphisms between the objects
of $\HC_\diag$, in terms of light leaf morphisms
\cite{LLL}. The following theorem shows that $\HC_\diag$ does indeed  give a ``generators and relations
description'' of the Hecke category:

\begin{thm}[{{\cite[Theorem 10.3.1]{RW}}}] \label{thm:equiv}
We have an equivalence of
  graded monoidal categories:
\[
\HC_{\diag}^{\Bbbk,\Kar} \simto \HC_{\geom}^{\Bbbk}.
\]
\end{thm}

\begin{remark}
  Knowing a presentation of a group or algebra by generators and
  relations opens the possibility of defining representations by
  specifying the action of generators and verifying
  relations. Similarly, studying actions of monoidal categories is
  sometimes easier when one has a presentation. In principle, the
  above presentation should allow a detailed study of categories acted
  on by the Hecke category. For interesting recent classification
  results, see \cite{MazorchukMiemietz2,MacTub}. One drawback of the
  theory in its current state is that the above relations (though
  explicit) can be difficult to check in examples. The parallel theory
  of representations of categorified quantum groups is much better
  developed (see e.g. \cite{CR,Brundan}).
\end{remark}

\begin{remark}
  An important historical antecedent to Theorems \ref{thm:diag_cat}
  and \ref{thm:equiv} is the theory of Soergel bimodules. We have
  chosen not to discuss this topic, as there is already a substantial
  literature on this subject. The above generators and relations were first written down in the
  context of Soergel bimodules, and Soergel bimodules are used in the
  proof of Theorem \ref{thm:equiv}. We refer the interested reader to 
  the surveys \cite{RBourbaki,LGentle} or the papers 
 \cite{S90,SHC,SB}. 
\end{remark}

\subsection{The spherical and anti-spherical module}

In this section we introduce the spherical and anti-spherical modules
for the Hecke algebra, as well as their categorifications. 
 They
are useful for (at least) two reasons: they are ubiquitous in
applications to representation theory; and they often provide smaller
worlds in which interesting phenomena become more tractable.

Throughout this section we fix a subset $I \subset S$ and assume for simplicity that the standard
parabolic subgroup $W_I$ generated by $I$ is finite. We denote by
$w_I$ its longest element. Let $H_I$ denote
the parabolic subalgebra of $H$ generated by $\delta_s$ for $s \in
I$; it is canonically isomorphic to the Hecke algebra of
$W_I$. Consider the induced modules
\[
M_I  :=  H\otimes_{H_I} \triv_v \quad
\text{and} \quad
N_I  := H \otimes_{H_I} \sgn_v
\]
where $\triv_v$ (resp. $\sgn_v$) is the rank one $H_I$-module with
action given by $\delta_s
\mapsto v^{-1}$ (resp. $\delta_s \mapsto -v$). These modules 
are the \emph{spherical} and \emph{anti-spherical} modules
respectively. If $W^I$ denotes the set of minimal length representatives
for the cosets $W/W_I$ then $\{ \delta_x \otimes 1  \; | \; x \in W^I \}$ gives
a (\emph{standard}) basis for $M_I$ (resp. $N_I$), which we
denote by $\{ \mu_x \;  | \; x \in W^I \}$ (resp. $\{ \nu_x \;  | \;
x \in W^I \}$). We denote the canonical bases in $M_I$ (resp. $N_I$)
by $\{ c_x \; | \; x \in W^I \}$ (resp. $\{ d_x \; | \; x \in W^I
\}$) (see e.g. \cite{SoeKL}).

We now describe a categorification of $M_I$. To $I$ is associated a
standard parabolic subgroup $\PS_I \subset \GS$, and we may
consider the partial flag variety $\GS/\PS_I$ (an ind-variety) and its Bruhat decomposition
\[
\GS/\PS_I  = \bigsqcup_{x \in W^I} Y_x \quad \text{where}
\quad Y_x :=\BS   \cdot x \PS_I /\PS_I .
\]
The closures $\overline{Y}_x$ are Schubert varieties, and we denote by
$\ic_{x,I}$ (resp. $\ESS_{x,I}$) the intersection cohomology complex
(resp. parity sheaf) supported on $\overline{Y}_x$.

Given any $\GS$-variety or ind-variety $Z$ the monoidal category
$D^b_\BS(\GS/\BS;\Bbbk)$ acts on $D^b_{\BS}(Z;\Bbbk)$. (The
definition is analogous to the formula for convolution given earlier.) In particular, $\HC_\geom^\Bbbk$ acts on
$D^b_{\BS}(\GS/\PS_I;\Bbbk)$. One can check that this action preserves
\[
\MC_I^{\Bbbk} := \langle \ESS_{x,I} \; | \; x \in W^I \rangle_{\oplus, [1]}
\]
and thus $\MC_I^{\Bbbk}$ is a module over $\HC_\geom^\Bbbk$. We have:

\begin{thm} There is a unique isomorphism of $H =
  [\HC_\geom^\Bbbk]_\oplus$-modules 
\[
M_I \simto [\MC_I^{\Bbbk}]_{\oplus}
\]
sending $\mu_{\id} \mapsto [ \Bbbk_{\PS_I /\PS_I}]$ (we use the indentification  $H =
  [\HC_\geom^\Bbbk]_\oplus$ of Theorem \ref{thm:geom_cat}).
\end{thm}

The inverse to the isomorphism in the theorem is given by the \emph{character map}
\begin{gather*}
  \ch : [\MC_I^{\Bbbk}]_{\oplus} \simto M_I \qquad
\FSS  \mapsto \sum_{x \in W^I} \dim_\ZM ( H^*(\FSS_{x\PS_I / \PS_I }) )
v^{-\ell(x)} \mu_x \in M_I.
\end{gather*}
(The notation is entirely analogous to the previous definition of
$\ch$ in \S\ref{hecke:geom}.)

We now turn to categorifying $N_I$. The full additive subcategory
\[
 \langle \ESS_x \; | \; x \notin W^I \rangle \subset \HC_\geom^\Bbbk
\]
is a left ideal. In particular, if we consider the quotient of
additive categories
\[
\NC_I^\Bbbk := \HC_\geom^\Bbbk / \langle \ESS_x \; | \; x \notin W^I \rangle
\]
this is a left $\HC^\Bbbk$-module. We denote the image of $\FSS \in
\HC_\geom^\Bbbk$ by $\overline{\FSS}_I$. The
objects $\overline{\ESS}_{x,I}$ for $x \in W^I$ are precisely the indecomposable
objects of $\NC^\Bbbk_I$ up to shift and isomorphism. We have:

\begin{thm} There is a unique isomorphism of right $H =
  [\HC_\geom^\Bbbk]_\oplus$-modules 
\[
N_I \simto [\NC_I^{\Bbbk}]_{\oplus}
\]
sending $\nu_{\id} \mapsto [ \overline{\ESS}_{\id,I}]$ (we use the identification  $H =
  [\HC_\geom^\Bbbk]_\oplus$ of Theorem \ref{thm:geom_cat}).
\end{thm}

The inverse $\ch :  [\NC_I^{\Bbbk}]_{\oplus} \simto N^I$ is more
complicated to describe.

\begin{remark}
  It is also possible to give a geometric description of $
  \NC_I^{\Bbbk}$ via Iwahori-Whittaker sheaves \cite[Chapter 11]{RW}.
\end{remark}

\subsection{The $p$-canonical basis} Suppose that $\Bbbk$ is a field
of characteristic $p \ge 0$. Consider the Hecke category
$\HC_\geom^\Bbbk$ with coefficients in $\Bbbk$. Let us define
\[
\p b_x := \ch(\ESS_x) \in H.
\]
Because $\ESS_x$ is supported on $\overline{X}_x$ and its restriction
to $X_x$ is $\Bbbk_{X_x}[\ell(x)]$, it follows from the
definition of the character map that 
\begin{equation}
  \label{eq:3}
  \p b_x = \delta_x + \sum_{y < x} \p h_{y,x} \delta_y
\end{equation}
for certain $\p h_{y,x} \in \ZM_{\ge 0}[v^{\pm 1}]$. Thus the set
$\{ \p b_x \; | \; x \in W \}$ is a basis for $H$, the
\emph{$p$-canonical basis}. The base change coefficients $\p h_{y,x}$
are called \emph{$p$-Kazhdan-Lusztig polynomials}, although they are Laurent
polynomials in general.

The $p$-canonical basis has the following properties (see \cite[Proposition 4.2]{JensenW}):
\begin{gather}
  \label{eq:self-dual}
d(\p b_x) = b_x \quad \text{for all $x
    \in W$;} \\ \label{eq:adj}
\text{if } \;\; \p b_x = \sum_{y \le x} \p a_{y,x}
    b_y \;\;\text{then $\p
    a_{y,x} \in \ZM_{\ge 0}[v^{\pm 1}]$ and $d(\p a_{y,x})=
    \p a_{y,x}$;} \\ \label{eq:mult}
\text{if } \;\; \p b_x \p b_y = \sum_{z \in W} \p \mu_{x,y}^z
    \p b_z \;\;\text{then $\p \mu_{x,y}^z\in \ZM_{\ge 0}[v^{\pm 1}]$ and $d(\p \mu_{x,y}^z) =
    \p \mu_{x,y}^z$;} \\
  \label{eq:largep}
\text{for fixed $x \in W$ we have $\p b_x = {}^0 b_x = b_x$ for large $p$.}
\end{gather}
\begin{remark} 
  Recall that the Kazhdan-Lusztig basis is uniquely determined by the
  ``self-duality'' and ``degree bound'' conditions (see
  \S\ref{sec:heckealg}). The $p$-canonical basis satisfies
  self-duality \eqref{eq:self-dual}, but there appears to be no
  analogue of the degree bound condition in general (see Example
  \ref{ex:degbound}   below).
\end{remark}

\begin{remark}
  There is an algorithm to calculate the $p$-canonical basis,
  involving the generators and relations presentation of the Hecke
  category discussed earlier. This algorithm is described in detail in
  \cite[\S 3]{JensenW}.
\end{remark}

\begin{remark}
  The Kazhdan-Lusztig basis only depends  on the Weyl
  group (a fact which is rather surprising from a geometric point of
  view). The $p$-canonical basis depends on the root 
  system. For example, the $2$-canonical bases in types $B_3$ and
  $C_3$ are quite different (see \cite[\S 5.4]{JensenW}).
\end{remark}

\begin{remark}
  The above properties certainly do not characterise the $p$-canonical
  basis. (For example, for affine Weyl groups the $p$-canonical bases
  are distinct for every prime.) However in certain
  situations they do appear to constrain the situation quite
  rigidly. For example, the above conditions are enough to deduce that
  $\p b_x = b_x$ for all primes $p$, if $\GS$ is of types $A_n$ for $n
  < 7$ (see \cite{W}). See \cite{ThorgeThesis}
  for further combinatorial constraints on the $p$-canonical basis.
\end{remark}

We can also define $p$-canonical bases in the spherical and
anti-spherical module.
Let $I \subset S$ be as in the previous section, and $\Bbbk$ and $p$
be as above. For $x \in W^I$, set
\begin{gather*}
\p c_x := \ch( \ESS_{x,I}) = \sum_{y \in W^I} \p m_{y,x} \mu_y \in M_I,  \\
\p d_x := \ch( \overline{\ESS}_{x,I}) = \sum_{y \in W^I} \p n_{y,x} \nu_y \in N_I.
\end{gather*}
We have $\p m_{y,x} = \p n_{y, x} = 0$ unless $y \le x$ and $\p m_{x,x} =
\p n_{x,x} = 1$. Thus $\{ \p c_x \}$ (resp. $\{ \p d_x \}$) give
\emph{$p$-canonical} bases for $M_I$ (resp. $N_I$). We leave it to
the reader to write down the
analogues of \eqref{eq:self-dual}, \eqref{eq:adj}, \eqref{eq:mult} and
\eqref{eq:largep} that they satisfy.

We define spherical and anti-spherical analogues of the ``adjustment'' polynomials $\p a_{y,x}$ via:
\begin{gather*}
  \p c_x = \sum_{y \in W^I}  \p a_{y,x}^\sph c_y \quad \text{and} \quad
  \p d_x = \sum_{y \in W^I}  \p a_{y,x}^\asph d_y.
\end{gather*}
These polynomials give partial information on the $p$-canonical
basis. For all $x, y \in W^I$ we have:
\begin{gather} \label{eq:sphasphrelations}
\p a_{yw_I, xw_I} = \p a^\sph_{y,x} \quad \text{and} \quad
\p a_{y, x} = \p a^\asph_{y,x}.
\end{gather}

We finish this section with a few examples of the $p$-canonical
basis. These are intended to complement the calculations in \S \ref{ex}.

\begin{ex} Let $\GS$ be of type $B_2$ with Dynkin diagram:
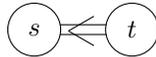
\[\begin{tikzpicture}[auto, baseline=(current  bounding  box.center), scale=1.3]
      \draw (0,\edgeShift) -- (-1,\edgeShift);
      \draw (0,-\edgeShift) -- (-1,-\edgeShift);
      \path (0,0) to node[Greater] (mid) {} (-1,0);
      \draw (mid.center) to +(30:\wingLen);
      \draw (mid.center) to +(330:\wingLen);
      \node [DynNode] (1) at (-1,0) {$s$};
      \node [DynNode] (2) at (0,0) {$t$}; 
\end{tikzpicture} \]
The Schubert variety $\overline{Y}_{st} \subset \GS/\PS_s$ has an
isolated singularity at $\PS_s /\PS_s$, and a neighbourhood of this
singularity is isomorphic to $X$ from Example \ref{ex:nullcone}. From
this one may deduce that
\[
{}^2 c_{st} = c_{st} + c_{\id}.
\]
For a version of this calculation using diagrams see \cite[\S
5.1]{JensenW}.
\end{ex}

\begin{ex} \label{ex:degbound}Here we explain the implications of
  Example \ref{ex:non-perverse} for the $p$-canonical basis.
The singularity $\CM^{2n}/(\pm 1)$ occurs in
the affine Grassmannian for $\Sp_{2n}$, which is isomorphic to
$\GS/\PS_I$, where $\GS$ is the affine Kac-Moody group 
of affine type $C_n$ with Dynkin diagram
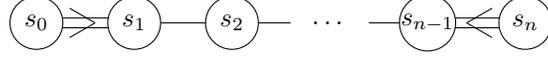
\[\begin{tikzpicture}[auto, baseline=(current  bounding  box.center),scale=1.3]
      \draw (0,\edgeShift) -- (-1,\edgeShift);
      \draw (0,-\edgeShift) -- (-1,-\edgeShift);
      \draw (3,\edgeShift) -- (4,\edgeShift);
      \draw (3,-\edgeShift) -- (4,-\edgeShift);
      \path (0,0) to node[Greater] (mid) {} (-1,0);
     \draw (-.4,0) to +(150:\wingLen);
      \draw (-.4,0) to +(210:\wingLen);
     \draw (3.4,0) to +(30:\wingLen);
      \draw (3.4,0) to +(330:\wingLen);
      \draw (0,0) -- (1,0);
      \draw (1,0) -- (1.6,0);
      \draw (2.4,0) -- (3,0);
      \node [DynNode] (1) at (-1,0) {$s_0$};
      \node [DynNode] (2) at (0,0) {$s_1$}; 
      \node [DynNode] (3) at (1,0) {$s_2$}; 
      \node (4) at (2,0) {$\dots$}; 
      \node [DynNode] (5) at (3,0) {${\tiny s_{n-1}}$}; 
      \node [DynNode] (6) at (4,0) {$s_n$}; 
\end{tikzpicture} \]
and $I = \{ s_1, \dots, s_n\}$ denotes the subset of finite reflections.
After some work matching parameters, one may deduce that
\[
{}^2 c_{w_n w_{n-1} s_0} = c_{w_n w_{n-1} s_0} + (v^{2n-2} + v^{2n-4} + \dots + v^{-2n+2})
\cdot c_{\id}.
\]
where $w_n$ (resp. $w_{n-1}$) denotes the longest element in the
standard parabolic subgroup generated by $\{ s_1, \dots, s_n \}$
(resp. $\{ s_2, \dots, s_n \}$).
\end{ex}

\begin{ex} Let $\GS = \SL_8(\CM)$ with simple reflections:
\[\begin{tikzpicture}[auto, baseline=(current  bounding  box.center),scale=1.3]
      \draw (0,0) -- (-1,0);
      \draw (0,0) -- (1,0);
      \draw (1,0) -- (2,0);
      \draw (2,0) -- (3,0);
      \draw (3,0) -- (4,0);
      \draw (4,0) -- (5,0);
      \node [DynNode] (1) at (-1,0) {$s_1$};
      \node [DynNode] (2) at (0,0) {$s_2$}; 
      \node [DynNode] (3) at (1,0) {$s_3$}; 
      \node [DynNode] (3) at (2,0) {$s_4$}; 
      \node [DynNode] (3) at (3,0) {$s_5$}; 
      \node [DynNode] (3) at (4,0) {$s_6$}; 
      \node [DynNode] (3) at (5,0) {$s_7$}; 
\end{tikzpicture} \]
Let
\[w = s_1s_3s_2s_4s_3s_5s_4s_3s_2s_1s_6s_7s_6s_5s_4s_3\]
and consider $w_I$ where $I = \{s_1, s_3,s_4,s_5,s_7 \}$. The singularity of
the Schubert variety $\overline{X}_w$ at $w_I$ is isomorphic to the
Kashiwara-Saito singularity from Example \ref{ex:KS} (with $d = 2$). It
follows that 
\[
{}^2 b_w = b_w + b_{w_I}.
\]
This is one of the first examples for $\SL_n$ with $\p b_x \ne b_x$.
\end{ex}

\subsection{Torsion explosion} \label{sec:explosion}
 In this section we assume that $\GS
\cong \SL_n(\CM)$ and so $W = S_n$, the symmetric group. Here the
$p$-canonical basis is completely known for $n = 2,3,\dots,9$ and
difficult to calculate beyond that. The following theorem makes clear
some of the difficulties that await us in high rank:

\begin{thm}[{{\cite{WT}}}] \label{thm:explosion}
Let $\un{\gamma}$ be a word of length $l$ in the generators
\[
\left ( \begin{matrix} 1 & 1 \\ 0 & 1 \end{matrix} \right ) \quad
\text{and} \quad
\left ( \begin{matrix} 1 & 0 \\ 1 & 1 \end{matrix} \right )
\]
with product
\[\gamma = \left ( \begin{matrix} \g_{11} & \g_{12} \\ \g_{21} &
    \g_{22} \end{matrix} \right ).\]
For non-zero $m \in \{ \g_{11},\g_{12},\g_{21}, \g_{22} \}$ and any
prime $p$ dividing $m$ there exists $y \in S_{3l + 5}$ such that
$\p b_y \ne b_y.$
\end{thm}

The moral seems to be that arithmetical issues (``which primes divide
entries of this product of elementary matrices?'') are hidden in the
question of determining the $p$-canonical basis.\footnote{Another
  example of this phenomenon from \cite{WT}: for any prime number $p$ dividing the $l^{th}$
  Fibonacci number there exists $y \in S_{3l + 5}$ with $\p b_y
  \ne b_y$. Understanding the behaviour of primes dividing Fibonacci numbers is
  a challenging open problem in number theory. It is conjectured, but not
  known, that infinitely many Fibonacci numbers are prime.}

We can get some qualitative information out of Theorem
\ref{thm:explosion} as follows. Define
\begin{gather*}
  \label{eq:1}
\Pi_n := \{ p \text{ prime }|\; \p b_x \ne b_x \text{ for some $x \in
  S_n$} \}.
\end{gather*}
Because any Schubert variety in $\SL_n(\CM)/\BS$ is also a Schubert
variety in $\SL_{n+1}(\CM)/\BS$ we have inclusions $\Pi_n \subset
\Pi_{n+1}$ for all $n$. By long calculations by Braden, Polo, Saito
and the author, we know the following about $\Pi_n$ for small $n$:
\begin{gather*}
\Pi_n = \emptyset \quad \text{for $n \le 7$,} \\
\Pi_n = \{ 2 \} \quad \text{for $n = 8,9$,} \\
\{ 2, 3\} \subset \Pi_{12}.
\end{gather*}
The most interesting values here are $2 \in \Pi_8$ (discovered by
Braden in 2002, see \cite[Appendix]{W}) and $3 \in \Pi_{12}$ (discovered by Polo in
2012). More generally, Polo shows that $p \in \Pi_{4p}$ for any prime
$p$, and hence $\Pi_n$ exhausts all prime numbers as $n \to \infty$
(see Example \ref{ex:KS}).

For applications to representation theory, it is important to know how
large the entries of $\Pi_n$ grow with $n$.\footnote{For example, the
  Lusztig conjecture would have implied that the entries of $\Pi_n$ are bounded
  linearly in $n$, and the James conjecture would have implied a
  quadratic bound in $n$, see \cite{WT}.} Some number theory,
combined with Theorem \ref{thm:explosion}, implies the following:

\begin{cor}[{{\cite[Theorem A.1]{WT}}}]
For $n$ large, $n \mapsto \max \Pi_n$ grows at least
  exponentially in $n$. More generally, $\Pi_n$ contains many
  exponentially large prime numbers.
\end{cor}

\begin{remark}
  Let us try to outline how Theorem \ref{thm:explosion} is
  proved. To $\un{\g}$ and $m$ we associate a reduced expression
  $\un{x} = (s_1, \dots, s_n)$ for some particular $x \in S_{3l + 5}$. (There is a
  precise but complicated combinatorial recipe as to how to do this,
  which we won't go into here. Let us mention however that the length
  of $x$ grows quadratically in $l$.) Associated to this
  reduced expression we have a Bott-Samelson resolution
\[
f : BS_{\un{x}} \to \overline{X}_x.
\]
We calculate the intersection form at a point $w_I\BS / \BS$
(corresponding to the maximal element of a standard parabolic
subgroup) and discover the $1 \times 1$-matrix $(m)$. Thus for any $p$
dividing $m$ the Decomposition Theorem fails for $f$ at the point
$w_I\BS / \BS$, which is enough to deduce the theorem. The
hard part in all of this is finding the appropriate expression
$\un{x}$ and calculating the intersection form. The intersection form calculation
was first done in \cite{WT} using a formula in the nil Hecke ring
discovered with He \cite{HeW}. Later a purely geometric argument was
found \cite{WIH}.
\end{remark}

\begin{remark}
  Let us keep the notation of the previous remark. In general we do
  not know whether $\p a_{w_A,x} \ne 0$ for any $p$ dividing $m$, only
  that there is some $y$ with $w_A \le y \le x$ and $\p a_{w_A,y}
  \ne 0$. Thus, in the statement
  of Theorem \ref{thm:explosion} we don't know that $\p b_x \ne
  b_x$, although this seems likely.
\end{remark}

\begin{remark}
  By a classical theorem of Zelevinsky \cite{Zelevinski}, Schubert
varieties in Grassmannians admit small resolutions, and hence the 
$p$-canonical basis is equal to the canonical basis in the spherical
modules for one step flag varieties (we saw a hint of this in Example
\ref{ex:gras}). It is an interesting question (suggested by Joe
Chuang) as to how
the $p$-canonical basis behaves in flag varieties with small numbers of
steps and at what point (i.e. at how many steps) the behaviour indicated in Theorem
\ref{thm:explosion} begins.
\end{remark}

\begin{remark}
  Any Schubert variety in $\SL_n(\CM)/\BS$ is isomorphic to a Schubert
  variety in the flag varieties of types $B_n$, $C_n$ and $D_n$. In
  particular, the above complexity is present in the $p$-canonical
  bases for all classical finite types.
\end{remark}

\subsection{Open questions about the $p$-canonical basis} In this
section we discuss some interesting open problems about the
$p$-canonical basis. We also try to outline what is known and point
out connections to problems in modular representation theory.

In the following a Kac-Moody root datum is assumed to be fixed
throughout. Thus, when we write $\p b_x$, its dependence on the root
datum is implicit. Throughout, $p$ denotes the characteristic of
$\Bbbk$, our field of coefficients.

\begin{question} \label{Q:one}
  For $x \in W$ and $p$ a prime, when is $\p b_x = b_x$?
\end{question}

\begin{remark}
  This question is equivalent to asking whether $\ic^\Bbbk_x
  \cong \ESS^\Bbbk_x$.
\end{remark}

A finer-grained version of this question is:

\begin{question} \label{Q:two}
  For $x,y \in W$ and $p$ a prime, when is $\p h_{y,x} = h_{y,x}$?
\end{question}

\begin{remark}
If $h_{y,x}  = v^{\ell(x) - \ell(y)}$ then Question \ref{Q:two} has a
satisfactory answer. In this case $y\BC/\BC$ is a rationally smooth
point of the Schubert variety $\overline{X}_x$ and $\p h_{y,x} =
h_{y,x}$ if and only if $\overline{X}_x$ is also $p$-smooth at
$y\BC/\BC$; moreover, this holds if and only if a certain
combinatorially defined integer (the numerator
of the ``equivariant multiplicity'') is not divisible by $p$, see
\cite{JWKumar, Dyer}. (See \cite{F3,FW} for related ideas.) It would be very interesting if one could
extend
such a criterion beyond the rationally smooth case.
\end{remark}

In applications the following variants of Question \ref{Q:one}
and \ref{Q:two} (for particular choices of $Z$) are more relevant:

\begin{question} \label{Q:ideal}
  Fix $Z \subset W$. For which $p$ does there exist $x \in Z$
  with $\p b_x \ne b_x$?
\end{question}

\begin{question} \label{Q:idealc}
  Fix $Z \subset W^I$.
  \begin{enumerate}
  \item For which $p$ is $\p c_x = c_x$ for all $x \in Z$?
  \item For which $p$ is $\p d_x = d_x$ for all $x \in Z$?
  \end{enumerate}
\end{question}

\begin{remark} \label{rem:LC}
  If $\GS$ is finite-dimensional then $\p b_x = b_x$ for all $x \in W$
  if and only if a part of Lusztig's conjecture holds (see
  \cite{Soe}). The results of \S\ref{sec:explosion} give
  exponentially large counter-examples fo the expected bounds in
  Lusztig's character formula \cite{L80,WT}.
\end{remark}

\begin{remark}
With $\GS$ as in the previous remark, Xuhua He has suggested that we
might have $\p b_x = b_x$ for all $x$, if $p > |W|$. This seems like a
reasonable hope, and it would be wonderful to have a proof.
\end{remark}

\begin{remark} \label{rem:LC}
  Suppose $\GS$ is an affine Kac-Moody group and $I \subset S$ denotes the
  ``finite'' reflections (so that $W_I = \langle I \rangle$ is the
  finite Weyl group). Then there exists a finite subset $Z_1 \subset W^I$
  for which Question \ref{Q:ideal}(1) is equivalent to determining
  in which characteristics Lusztig's character formula holds, see \cite[\S
  11.6]{ARFM} and \cite[\S 2.6]{W-Takagi}. Because $Z_1$ is finite,
  $\p c_x = c_x$  for all $x \in Z_1$ for $p$ large, which
  translates into the known fact that Lusztig's conjecture holds in
  large characteristic.
\end{remark}

\begin{remark}
  Suppose that $\GS$ and $I$ are as in the previous remark.
  There exists a subset $Z(p) \subset  W^I$ (depending
  on $p$) such that if $\p d_x = d_x$ for all $x \in Z(p)$ then
  Andersen's conjecture on characters of tilting modules (see
  \cite[Proposition 4.6]{AFilt2})  holds in characteristic
  $p$ (this is a consequence of the character formula proved in
  \cite{AMRW2}). Note that Andersen's conjecture does not give a character
  formula for the characters of all tilting modules and is not known even
  for large $p$.
\end{remark}

The following questions are also interesting:

\begin{question} \label{Q:perverse}
  For which $x \in W$ and $p$ is $\p a_{y,x} \in \ZM$ for all $y \in W$?
\end{question}

\begin{remark}
This is equivalent to asking when $\ESS_x^\Bbbk$ is perverse.
\end{remark}

\begin{question} \label{Q:perverse2}
  Fix $? \in \{ \sph,\asph \}$. For which $x \in W^I$ and $p$ is $\p a^?_{y,x}
  \in \ZM$ for all
  $y \in W^I$?
\end{question}

\begin{remark}
  Example \ref{ex:degbound} shows that in the affine case $\p a_{y,x}$ can be a polynomial in
  $v$ of arbitrarily high degree.
An example of Libedinsky and the author
  \cite{LWnp} shows that there exists $x,y$ in the symmetric group
  $S_{15}$, for which ${}^2 a_{y,x} = (v+v^{-1})$. Recently
  P.~McNamara has proposed new candidate examples, which appear to
  show that for any $p$, the degree of $\p a_{y,x}$ is unbounded in
  symmetric groups.
\end{remark}

\begin{remark}
  Suppose $\GS$ and $I$ are as in Remark \ref{rem:LC}. It follows from
  the the results of \cite{JMW3,MR} that $\p a_{y,x} \in \ZM$ if $x$
  is maximal in $W_I x W_I$. More generally, it seems likely that $\p
  a_{y,x}^\sph \in \ZM$ for all $x, y \in W^I$ and large $p$ (depending only on the Dynkin
  diagram of $\GS$). This is true for trivial reasons in affine
  types $A_1$ and $A_2$.
\end{remark}

\begin{remark}
  In contrast, recent conjectures of Lusztig and the author
  \cite{LWbilliards}
  imply that, if $\GS$ is of affine type $\widetilde{A_2}$ then it is
  never the case (for any $p \ne 2$) that $\p a^\asph_{y,x} \in \ZM$ for all $y, x \in W^I$. In fact, our conjecture implies that
  \[
    \max \; \{ \deg( \p a^\asph_{y,x}) \; | \; y, w \in W^I \} = \infty.
  \]
This contrast in behaviour between the
  $p$-canonical bases in the spherical and anti-spherical modules is
  rather striking. 
\end{remark}

\section{Koszul duality} \label{sec:KD}

In this section we discuss Koszul duality for the Hecke category. This
is a remarkable derived equivalence relating the Hecke categories
of Langlands dual groups. It resembles a Fourier
transform. Its modular version involves parity sheaves,
and is closely related to certain formality 
questions. In this section we assume that the reader has some
background with perverse sheaves and highest weight categories.

\subsection{Classical Koszul duality} 
Let $C, \GS, \BS, \TS, W,\Bbbk$ be as previously. 
We denote by $\GS^\vee, \BS^\vee, \TS^\vee$ the Kac-Moody group
(resp. Borel subgroup, resp. maximal torus) associated to the
 dual Kac-Moody root datum. We have a canonical identification of $W$
 with the Weyl group of $\GS^\vee$.

In this section we assume that $\GS$ is a (finite-dimensional)
complex reductive group, i.e. that $C$ is a Cartan matrix. We denote by $w_0 \in W$ the longest
element. For any
$x \in W$ let $i_x : X_x = \BS \cdot x\BS/\BS \into \GS/\BS$ denote the inclusion of
the Schubert cell and set
\[
\Delta_x := i_{x!} \Bbbk_{X_x}[\ell(x)] \quad \text{and} \quad
\nabla_x := i_{x*} \Bbbk_{X_x} [\ell(x)].
\]
Let $D^b_{(\BS)}(\GS/\BS;\Bbbk)$ denote the derived category,
constructible with respect to $\BC$-orbits and let
$\Perv_{(\BS)}(\GS/\BS;\Bbbk) \subset D^b_{(\BS)}(\GS/\BS;\Bbbk)$
denote the
subcategory of perverse sheaves. The abelian category
$\Perv_{(\BS)}(\GS/\BS;\Bbbk)$ is highest weight \cite{BGS,BBM} with standard
(resp. costandard) objects $\{ \Delta_x \}_{x \in W}$ (resp. $\{
\nabla_x \}_{x \in W}$). For $x \in W$, we denote  by $\PSS_x, \ISS_x$ and
$\TSS_x$ the corresponding indecomposable projective, injective,
and tilting object. The corresponding objects in 
$\Perv_{(\BS^\vee)}(\GS^\vee/\BS^\vee;\Bbbk)$ are denoted with a check,
e.g. $\ic_x^\vee, \Delta_x^\vee$ etc.


Let us assume that $\Bbbk = \QM$. Motivic considerations, together with
the Kazhdan-Lusztig inversion formula (see \cite{KLp})
\begin{equation}
  \label{eq:KLinversion}
\sum_{z \in W} (-1)^{\ell(x) + \ell(y)}h_{y,x}h_{yw_0,zw_0} = \delta_{x,z},
\end{equation}
led Beilinson and Ginzburg \cite{BGp} to the following conjecture\footnote{to
  simplify the exposition we have
  modified the statement of their original conjecture slightly (they
  worked with Lie algebra representations and sought a contravariant
  equivalence).}:
\begin{enumerate}
\item There exists a triangulated category $D^{mix}_{(\BS)}(\GS/\BS; \QM)$
  equipped with an action of the integers $\FSS \mapsto \FSS\langle m
  \rangle$ for $m \in \ZM$ (``Tate twist'') and a ``forgetting the mixed structure'' functor
\[
\phi :
  D^{mix}_{(\BS)}(\GS/\BS;\QM) \to D^b_{(\BS)}(\GS/\BS;\QM),\]
such that
\[
\Hom(\phi(\FSS),\phi(\GSS)) = \bigoplus_{n \in \ZM}  \Hom(\FSS,\GSS \langle
n \rangle)
\]
for all $\FSS, \GSS \in
D^{mix}_{(\BS)}(\GS/\BS;\QM)$. Furthermore, ``canonical'' objects
(e.g. simple, standard, projective etc.~objects) admit
lifts\footnote{$\FSS \in D^b_{(\BS)}(\GS/\BS;\QM)$ admits a lift, if
  there
  exists $\widetilde{\FSS} \in D^{mix}_{(\BS)}(\GS/\BS;\QM)$ such that $\FSS
  \cong \phi(\widetilde{\FSS})$.} to $D^{mix}_{(\BS)}(\GS/\BS;
\QM)$.
\item There is an equivalence of
triangulated categories
\begin{equation}
  \label{eq:BGS}
\kappa : D^{mix}_{(\BS)}(\GS/\BS; \QM) \simto D^{mix}_{(\BS^\vee)}(\GS^\vee/\BS^\vee; \QM)  
\end{equation}
such that $\kappa \circ \langle -1 \rangle[1] \cong \la 1 \ra \circ
\kappa$, and such that $\kappa$ acts on standard, simple and
projective objects (for an appropriate choice of lift) as follows:
\[
\Delta_x \mapsto \nabla^\vee_{x^{-1}w_0}, \qquad
\ic_x \mapsto \ISS^\vee_{x^{-1}w_0}, \qquad
\PSS_x \mapsto \ic^\vee_{x^{-1}w_0}.
\]
\end{enumerate}

\begin{remark}
  To understand why the extra grading (provided by the mixed structure)
  as well as the relation $\kappa \circ \langle -1 \rangle[1] \cong
  \la 1 \ra \circ \kappa$ is necessary, one only needs to ask oneself
where the grading on extensions between simple modules
should go under this equivalence.
\end{remark}

\begin{remark}
One can deduce from \eqref{eq:KLinversion} and the Kazhdan-Lusztig conjecture
  that the assignment $\Delta_x \mapsto \nabla_{x^{-1}w_0}$ on mixed
  categories forces $\ic_x \mapsto \ISS_{x^{-1}w_0}$ and $\PSS_x \mapsto
  \ic_{x^{-1}w_0}$ on the level of Grothendieck groups.
\end{remark}

This conjecture was proved by Beilinson, Ginzburg and Soergel in the
seminal paper \cite{BGS}, where they interpreted
$\kappa$ in the framework of Koszul duality for graded algebras. The
authors give two constructions of the mixed derived
category: one involving mixed \'etale sheaves (here it is necessary to
consider the flag variety for the split group defined over a finite field), and one involving mixed Hodge
modules.

\begin{remark}  Both constructions of the mixed
  derived category in \cite{BGS} involve some non-geometric
  ``cooking'' to get the right result.  Recently Soergel and Wendt have used various flavours of
  mixed Tate motives to give a purely geometric construction of these
  mixed derived categories \cite{SW16}.
\end{remark}

After the fact, it is not difficult to see that the mixed derived
category admits a simple definition. Indeed, the results of \cite{BGS}
imply that one has an equivalence
\begin{equation}
  \label{eq:SS}
\sigma : K^b(\ses(\GS/\BS;\QM)) \simto D^{mix}_{(\BS)}(\GS/\BS; \QM) .
\end{equation}
Here $\ses(\GS/\BS;\QM)$ denotes the full additive subcategory of
$D^{mix}_{(\BS)}(\GS/\BS; \QM)$ consisting of direct sums of shifts of
intersection cohomology complexes (``semi-simple complexes''), and
$K^b(\ses(\GS/\BS;\QM))$ denotes its homotopy category. Note that there are
two shift functors on $\ses(\GS/\BS;\QM)$: one coming from its structure as
a homotopy category (which we denote $[1]$); and one induced from the
shift functor on $\ses(\GS/\BS;\QM)$ (which we rename $(1)$).  Under the
equivalence $\sigma$, Tate twist $\langle 1 \rangle$ corresponds to
$[1](-1)$.

Now, if $\HC^\QM$ denotes the Hecke category we have an equivalence
\[
\QM \otimes_R \HC_\diag^\QM \simto \ses(\GS/\BS;\QM).
\]
Moreover, the left hand side can be described by generators and
relations. In particular, Koszul duality can be formulated
entirely algebraically as an equivalence
\[
\kappa: K^b(\QM \otimes_R \HC_\diag^\QM) \simto K^b(\QM \otimes_R \HC_\diag^{\vee,\QM}).
\]
The existence of such an equivalence (valid more generally for any finite
real reflection group, with $\QM$ replaced by $\RM$) has recently been
established by Makisumi 
\cite{MakiK}. (The case of a
dihedral group was worked out by Sauerwein \cite{Sauerwein}.)

\subsection{Monoidal Koszul duality} 
The above results raise the following questions:
\begin{enumerate}
\item How does Koszul duality interact with the monoidal structure?
\item Can Koszul duality be generalised to the setting of Kac-Moody
  groups?
\end{enumerate}
The first question was addressed by Beilinson and Ginzburg
\cite{BGwallcrossing}. They noticed that if one composes Koszul
duality $\kappa$ with the Radon transform and inversion, one obtains
a derived equivalence
\begin{equation}
  \label{eq:BG}
\widetilde{\kappa} : D^{mix}_{(\BS)}(\GS/\BS; \QM) \simto D^{mix}_{(\BS^\vee)} (\BS^\vee
\setminus \GS^\vee; \QM)  
\end{equation}
with $\widetilde{\kappa}  \circ \langle -1 \rangle[1] \cong \la 1 \ra \circ \widetilde{\kappa} $
as previously, however now 
\begin{equation}\label{eq:neweffect}
\ic_x \mapsto \TSS_x^\vee, \quad \Delta_x \mapsto \Delta_{x}^\vee, \qquad
\nabla_x \mapsto \nabla_{x}^\vee, \qquad
\TSS_x \mapsto \ic_{x}^\vee.
\end{equation}
The new equivalence $\widetilde{\kappa}$ is visibly more symmetric
than $\kappa$. It also has the advantage that it does not involve the
longest element $w_0$, and hence makes sense for Kac-Moody groups.
Moreover, Beilinson and Ginzburg conjectured that $\widetilde{\kappa}$ can be promoted to a
monoidal equivalence (suitably interpreted).

\begin{remark} It has been a stumbling block for some time
  that \eqref{eq:BG} cannot be upgraded to a monoidal equivalence in a
  straightforward way. This is already evident for $\SL_2$: the ``big'' tilting
  sheaf $\TSS_s \in D^b_{(\BS)}(\PM^1; \QM)$ does not admit a
  $\BS$-equivariant structure.
\end{remark}

Subsequently, Bezrukavnikov and Yun \cite{BY} established a monoidal equivalence
\begin{equation}
  \label{eq:BY}
\widetilde{\kappa} : (D^{mix}_\BS(\GS/\BS; \QM), *) \simto
(\widehat{D}^{mix}(\BS^\vee\leftdash \GS^\vee\rightdash \BS^\vee; \QM), \star)
\end{equation}
which induces the Koszul duality equivalence above after killing the
deformations, and is valid for any Kac-Moody group.\footnote{Actually,
  Bezrukavnikov and Yun use mixed $\ell$-adic sheaves, and no
  non-geometric ``cooking'' is necessary.} 
Here $\widehat{D}^{mix}(\BS^\vee\leftdash \GS^\vee\rightdash \BS^\vee; \QM)$
denotes a suitable (``free monodromic'') completion of the full
subcategory of mixed $U^\vee$-constructible complexes on $\GS^\vee/U^\vee$
which have unipotent monodromy along the fibres of the map
$\GS^\vee/U^\vee \to \GS^\vee/\BS^\vee$. The construction of this completion
involves considerable technical difficulties. The proof involves
relating both sides to a suitable category of Soergel bimodules
(and thus is by ``generators and relations'').


\subsection{Modular Koszul duality}

We now discuss the question of how to generalise \eqref{eq:BG} to
coefficients $\Bbbk$ of positive characteristic.  A first difficulty
is how to make sense of the mixed derived
category. A naive attempt (carried out in \cite{RSW}) is to consider a flag
variety over a finite field together with the Frobenius endomorphism
and its weights, however here one runs into problems because one
obtains gradings by a finite cyclic group rather than $\ZM$. Achar
and Riche took the surprising step of simply defining
\[
D^{mix}_{(\BS)}(\GS/\BS;\Bbbk) := K^b(\Parity(\GS/\BS;\Bbbk))
\]
where $\Parity(\GS/\BS;\Bbbk)$ denotes the additive category of
$\BS$-constructible parity complexes on $\GS/\BS$, and the shift $[1]$
and twist $(1)$ functors are defined as in the paragraph
following \eqref{eq:SS}. (The discussion there shows that this definition is
consistent when $\Bbbk = \QM$.) In doing so one obtains a triangulated category with
most of the favourable properties one expects from the mixed derived
category. In this setting Koszul duality takes the form:

\begin{thm}
\label{thm:mkd}
There is an equivalence of triangulated categories
\[
\kappa : \Dmix_{(\BS)}(\GS/\BS;\Bbbk) \simto \Dmix_{(\BS^\vee)}(\BS^\vee \backslash \GS^\vee;\Bbbk)
\]
which satisfies $\kappa \circ \langle -1 \rangle[1] \cong \la 1 \ra \circ \kappa$ and
\[
\kappa(\Delta_w) \cong \Delta^\vee_w, \quad \kappa(\nabla_w) \cong \nabla^\vee_w, \quad \kappa(\ESS_w) \cong \TS^\vee_w, \quad \kappa(\TS_w) \cong \ESS^\vee_w.
\]
\end{thm}

\begin{remark}
  The important difference  in the modular case is that tilting sheaves correspond to parity sheaves (rather than IC sheaves).
\end{remark}

\begin{remark}
  For finite-dimensional $\GS$ this theorem was proved in
  \cite{ARMixed2} (in good characteristic).   For general $\GS$ this theorem is proved in
  \cite{AMRW1,AMRW2}, as a corollary of a monoidal modular Koszul
  duality equivalence, inspired by \cite{BY}.
\end{remark}

\begin{remark}
The appearance of the Langlands dual
group was missing from the
original conjectures of Beilinson-Ginzburg \cite{BGp} and only
appeared in \cite{BGS}. However in the settings considered there
($\Bbbk = \QM$), the
Hecke categories associated to dual groups are equivalent. This is no
longer the case with modular coefficients, and examples (e.g. $B_3$
and $C_3$ in
characteristic 2) show that the analogue of Theorem \ref{thm:mkd}  is
false if one ignores the dual group.
\end{remark}

\begin{remark}
  A major motivation for   \cite{AMRW1,AMRW2} was a conjecture of
  Riche and the author \cite[\S 1.4]{RW} giving characters for tilting modules for
  reductive algebraic groups in terms of $p$-Kazhdan-Lusztig
  polynomials. In fact, a recent theorem of Achar and Riche \cite{ARFM}
  (generalising a theorem of Arkhipov, Bezrukavnikov and
  Ginzburg \cite{ABG}) combined with (a variant of) the above Koszul
  duality theorem leads to a solution of this conjecture. We expect
  that modular Koszul duality will have other applications
  in modular representation theory.
\end{remark}

\begin{remark} One issue with the above definition of the mixed
  derived category is the absence of a ``forget the mixed structure'' functor 
$\phi : D^{mix}_{(\BS)}(\GS/\BS;\Bbbk) \to
D^b_{(\BS)}(\GS/\BS;\Bbbk)$ in general. For finite-dimensional $\GS$ its
existence is established in \cite{ARMixed2}. Its existence for affine Weyl
groups would imply an important conjecture of Finkelberg and
Mirkovi\'c \cite{FM, ARFM}.
\end{remark}


\def\cprime{$'$} \def\cprime{$'$} \def\cprime{$'$}


\end{document}